\documentclass[twoside]{amsart} \usepackage{amsmath}
\usepackage{amssymb}
\usepackage{rotating}
\usepackage[colorlinks=true]{hyperref}
\usepackage{lscape}
\newcommand{\Hom}{\mbox{\rm Hom}}
\newcommand{\Ind}{\mbox{\rm Ind}}
\newcommand{\Ccal}{\ensuremath{\mathcal{C}}}
\newcommand{\Bcal}{\ensuremath{\mathcal{B}}}
\newcommand{\Dcal}{\ensuremath{\mathcal{D}}}

\begin{document}

\begin{abstract}
We determine the $2$-modular and $3$-modular character tables of the 
sporadic simple Harada-Norton group and its automorphism group.
\end{abstract}


\title{The Brauer characters of the sporadic simple
Harada-Norton group 
and its automorphism group in characteristics~$2$ and~$3$}

\author{Gerhard Hiss, J\"urgen M\"uller, Felix Noeske, and Jon Thackray}

\dedicatory{Dedicated to the memory of Herbert Pahlings}

\address{Lehrstuhl D f\"ur Mathematik, RWTH Aachen University,
52056 Aachen, Germany}

\email{G.H.: Gerhard.Hiss@math.rwth-aachen.de}
\email{J.M.: Juergen.Mueller@math.rwth-aachen.de}
\email{F.N.: Felix.Noeske@math.rwth-aachen.de}
\email{J.T.: jgt@pobox.com}

\maketitle

\section{Introduction}
\markboth{MODULAR CHARACTERS OF THE HARADA-NORTON GROUP}
{GERHARD HISS, J\"URGEN M\"ULLER, FELIX NOESKE, AND JON THACKRAY}

The present paper is a contribution to the modular ATLAS project
\cite{web:MODATLAS}. 
The aim of the project is to classify the modular
representations 
of the simple groups and their bicyclic
extensions given in the ATLAS \cite{ATLAS}, by computing their 
Brauer characters,
or equivalently, their decomposition numbers.

In particular for the sporadic simple groups this is a challenging task.
In lieu of a unified theory for these groups, the solution of this
problem has always relied on computational methods, and especially on
the MeatAxe (see \cite{MeatAxe}) developed by Parker and Thackray in the
late 1970's. However, even though the MeatAxe has evolved into a host of
programmes to analyse modules of matrix algebras, the open problems in
the modular ATLAS project defy a direct application:
owing to the magnitude of the computational problems, memory and time
constraints render a naive approach utterly infeasible.

To regain computational tractability, Parker and Thackray established
the condensation method (see \cite{Thackray:PhD}). 
The idea of condensation is that given a module
$V$ for an algebra $A$, one considers the \emph{condensed} module $Ve$
of the algebra $eAe$ for some suitably chosen idempotent $e \in A$. The
condensation functor from 
\textsf{mod}-$A$ to \textsf{mod}-$eAe$ has a
number of interesting properties, among the most important of which is
that it maps simple modules to simple modules or zero, and 
composition series to 
series whose layers are simple or zero.
For a detailed introduction to condensation see
\cite{Ryba90}, for example. 
Unfortunately, the use of condensation introduces the \emph{generation
problem} (see \cite{tackling}, for example): 
for computational purposes we need to work with a feasibly
small generating set for the \emph{condensed algebra} $eAe$, yet such
a set is not known in general, and we therefore devote time and energy to
verify our condensation results.

\medskip
{\bf Results.}
In this paper we determine the decomposition numbers of the sporadic simple 
Harada-Norton group $\mathsf{HN}$ and its automorphism group $\mathsf{HN}.2$
in characteristics 2 and 3. 
More precisely, the decomposition matrices of the $2$-blocks of
positive defect of $\mathsf{HN}$ and $\mathsf{HN}.2$ are
determined in Section~\ref{sec:2}, where the results for $\mathsf{HN}$ 
are given in
Tables~\ref{DecMatB1mod2} and \ref{tab:decB0mod2}, 
and those for $\mathsf{HN}.2$ are given in 
Tables~\ref{tab:decmatB0HN2mod2} and \ref{tab:decmatB12HN2mod2}. For
the reader's convenience we give the degrees of Brauer characters lying in
the principal $2$-block of $\mathsf{HN}$ in Table~\ref{tab:degsB0mod2}.
\begin{table}
 \centering
 \begin{tabular}{rrrrr}
  1& 132& 132& 760& 2\,650\\ 
  2\,650& 3\,344& 15\,904& 31\,086& 31\,086\\ 
  34\,352& 34\,352& 43\,416& 43\,416& 177\,286\\ 
  217\,130& 1\,556\,136&
 \end{tabular}
 \caption{Brauer character degrees in the principal $2$-block of
 \textsf{HN}}
 \label{tab:degsB0mod2}
\end{table}
The characteristic 3 case is dealt with in Section~\ref{sec:mod3},
where the results for $\mathsf{HN}$ are given in
Tables~\ref{DecMatB2}, \ref{tab:decB0mod3} and \ref{DecMatB1},
and those for the principal block of $\mathsf{HN}.2$, being
the only block not Morita equivalent to a block of $\mathsf{HN}$,
are given in Table~\ref{tab:decB0HN2mod3}. The degrees of the
irreducible Brauer characters in the principal $3$-block of $\mathsf{HN}$
are restated in Table~\ref{tab:degsB0mod3}.
\begin{table}
 \centering
 \begin{tabular}{rrrrr}
  1& 133& 133& 760& 3\,344\\ 
  8\,778& 8\,778& 9\,139& 12\,264& 12\,264\\
  31\,768& 31\,768& 137\,236& 147\,061& 255\,037\\
  339\,702& 339\,702& 496\,924& 496\,924& 783\,696
 \end{tabular}
 \caption{Brauer character degrees in the principal $3$-block of \textsf{HN}}
 \label{tab:degsB0mod3}
\end{table}

The results obtained here are also accessible in \cite{web:MODATLAS} 
and will be through the \textsf{CTblLib} database \cite{ctbllib} 
in \textsf{GAP} \cite{GAP4}.
Since the Brauer trees of $\mathsf{HN}$ and $\mathsf{HN}.2$ 
have already been determined in \cite{BrauerBaumBuch} and
the characteristic 5 case has been tackled in \cite{LuxNoeskeRyba},
this completes the modular ATLAS project 
for $\mathsf{HN}$. Moreover, since the decomposition matrix
of the non-principal $3$-block of $\mathsf{HN}$ of defect $2$ was
known prior to the completion of the present paper, it was the
starting point of further investigations into this block carried
out in \cite{KosMue}.

\medskip\noindent
The {\bf proofs} for the two characteristics are quite different:

In {\bf characteristic 2}, the non-principal block may be treated entirely
theoretically. It possesses a semi-dihedral
defect group of order 16, and we may therefore use Erdmann's
classification \cite{ErdmannSemidihedral} 
of possible decomposition matrices for such blocks
to determine its decomposition matrix. The
true challenge lies in the principal block. Here we use condensation,
but since there are no condensation subgroups readily available which
allow the application of \cite[Theorem 2.7]{tackling}, we are forced to
verify the condensation results. As the dimensions of the simple modules
with the exception of the last are small enough, we construct bases
for these simple modules in order to show that the degrees of the Brauer
atoms produced (see \cite[Definition 3.2.1]{manus:MOC}), 
which are lower bounds for the degrees of the
Brauer characters, are maximal, i.e.\ the atoms are in fact characters.
We can then show that the last atom is a character, too.

For the automorphism group the result is immediate from the ordinary 
character table.

In {\bf characteristic 3}, we start by investigating the block of defect
2 by choosing a basic set of projective characters and improving it
pursuing a strategy first employed in \cite{MueHabil}:
We consider the endomorphism ring $E$ of a suitably chosen 
projective lattice, and by Fitting correspondence 
compare the decomposition numbers of $E$ with those of the
block in question (see \cite[Ch.II.12]{LandrockBook}).
>From the practical side, we consider a direct summand
of a projective permutation module, so that the regular representation 
of $E$ is just the associated condensed module. 
Actually, here we are better
off with the generation problem: since the regular action of $E$
is faithful we are able to determine the dimension
of the condensed algebra in question. 
Alone, this strategy turns out to be  
computationally tractable only for projective modules
leading to non-faithful trace idempotents, 
so that after this analysis we are still left with 
several possible cases for the decomposition matrix. 

These cases are then greatly reduced in parallel to our treatment
of the principal block, which relies entirely on condensation. 
There are several condensation subgroups available which
are normal in maximal subgroups, making them ideal candidates for
\cite[Theorem 2.7]{tackling}. However, we still do not find any 
single suitable one whose trace idempotent is faithful.
Therefore, we choose two condensation subgroups
such that no simple module is annihilated by both condensation
idempotents. This introduces a new problem: given the
simple modules of one condensed algebra, we have to match them to the
simple modules of the other condensed algebra, in order to gain the
complete picture of the composition factors of the original module.
We use \cite{matching2} to overcome this problem (see also \cite{matching}).

For the automorphism group $\mathsf{HN}.2$, a simple application of Clifford
theory already gives the bulk of its Brauer characters. The remaining ones
are once more determined through condensation, reusing the
condensation infrastructure we have established for the simple group.

\medskip
{\bf Background.} 
In the sequel, we assume the reader to be familiar with the basics
of ordinary and modular group representation theory, and in
particular its computational aspects, as are for example exposed in
\cite{manus:MOC,LuxPah}. We employ the particular choice of Brauer
lifts using Conway polynomials as in \cite{bookModatlas}. Character
theoretic computations are carried out in \textsf{GAP}, whose facilities
are also used to explicitly compute with permutation and matrix
groups; in particular we use the \textsf{orb} package \cite{orb}.
Module theoretic computations are done using the MeatAxe and its
extensions \cite{CMeatAxe}, 
and as a computational workhorse we apply condensation
of tensor products, respectively homomorphism spaces, as described
in \cite{art:newtensorcondense} and \cite{LuxWie}. For a general
introduction to condensation we refer the reader to \cite{Ryba90}, and
to \cite{LuxMuellerRinge} for a treatise on peakword condensation.
As a source of explicit data we use the database \textsf{CTblLib} in
\textsf{GAP}, as well as the web-based database \cite{web:ATLAS},
also accessible via the \textsf{GAP} package \textsf{AtlasRep}
\cite{AtlasRep}. Next to matrix and permutation representations
for many of the groups in \cite{ATLAS}, given in \emph{standard
generators} in the sense of \cite{stdgens} where available, the latter
in particular also contains \emph{straight line programmes} to find
maximal subgroups and conjugacy class representatives.

\medskip
Let us fix some notation which will hold throughout the paper.
Let $G$ denote the Harada-Norton group $\mathsf{HN}$. Characters,
representations, and modules are given by their degrees, respectively
dimensions. By a slight abuse of notation, we therefore do not
distinguish between a module, the representation it affords, and the
associated Brauer character. If there is more than one character of
the same degree we affix subscripts to distinguish them. The ordinary
characters are numbered as they are numbered in the character tables
of \textsf{CTblLib}, which generally corresponds to the numbering
given in \cite{ATLAS}. In the decomposition matrices we also indicate
this numbering next to the character degrees if necessary. If $\chi$
is an ordinary character we let $\chi'$ denote its restriction to the
$p$-regular classes, where $p$ is the characteristic of the field
currently considered. Moreover, we stick to group theoretic notation as
is used in \cite{ATLAS}.

>From the ordinary character table of $G$ and $G.2$, see \cite{ATLAS},
we see that all character fields are $2$-elementary abelian,
thus by \cite[Section~I.5]{bookModatlas} the quadratic extension
$\mathrm{GF}(p^2)$ of any finite prime field $\mathrm{GF}(p)$ is a
splitting field for $G$ and $G.2$. We let $F$ denote $\mathrm{GF}(p^2)$
in the sequel, and let $(Q,R,F)$ denote a $p$-modular splitting system
for $G$ and $G.2$.

\section{Characteristic~$2$}
\label{sec:2}
The $2$-modular characters of $G$ are distributed over three blocks,
two of which have positive defect: the principal block $B_0$ of defect
$d=14$, having $k=45$ ordinary and $l=17$ modular characters, and the
block $B_1$ of defect $d=4$, having $k=8$ ordinary and $l=3$ modular
characters. The block of defect 0 consists of the ordinary character
$\chi_{46}$.

\subsection{Proof for the block $B_1$}
By \cite[Section 9]{Landrock} the defect group of Block
$B_1$ is a semidihedral group of order 16. In \cite[Section
11]{ErdmannSemidihedral} Erdmann lists all possible decomposition
matrices for such blocks, such that we are left with to decide which of
the four cases in Lemmas~11.4, 11.6, 11.9, and 11.11 holds:

In Erdmann's notation \cite[Theorem 11.1]{ErdmannSemidihedral},
work of Olsson \cite{olsson} implies that
$B_1$ has four characters $\chi_1,\ldots,\chi_4$ of height zero, 
three characters $\chi_1^*,\ldots,\chi_3^*$ of height one for which 
$\chi^{*\prime}:=\chi_1^{*\prime}=\chi_2^{*\prime}=\chi_3^{*\prime}$,
and a single character $\hat\chi$ of height two.
Moreover, $\chi_1,\ldots,\chi_4$ can be chosen such that 
$$ \begin{array}{rcrcl}
\delta_1\chi_1'+\delta_2\chi_2' &=& -\delta_3\chi_3'-\delta_4\chi_4' &=& 
\chi^{*\prime}, \\
\delta_2\chi_2'+\delta_4\chi_4' &=& -\delta_1\chi_1'-\delta_3\chi_3' &=& 
\hat\chi', \\
\end{array} $$
where $\delta_1,\ldots,\delta_4\in\{\pm 1\}$. 
In particular, choosing any $\chi_i\in\{\chi_1,\ldots,\chi_4\}$
and $\chi_r^*\in\{\chi_1^*,\ldots,\chi_3^*\}$, we conclude
that $\{\chi_i',\chi_r^{*\prime},\hat\chi'\}$ 
is a basic set of Brauer characters.
Actually, as is easily seen, the above conditions determine
$\chi_1,\ldots,\chi_4$ and $\delta_1,\ldots,\delta_4$ uniquely
up to the permutation interchanging
$\delta_1\chi_1\leftrightarrow -\delta_4\chi_4$ and
$\delta_2\chi_2\leftrightarrow -\delta_3\chi_3$,
and by the positivity of character degrees we are left with the four cases
$$ [\delta_1,\ldots,\delta_4]\in
\{[1,1,-1,-1],[1,-1,-1,1],[-1,1,-1,1],[1,1,-1,1]\},
$$
corresponding to the Lemmas cited above.

By computing the actual heights of the characters in $B_1$, we conclude that 
$\{\chi_1,\ldots,\chi_4\} = \{\chi_{17},\chi_{37},\chi_{45},\chi_{49}\}$,
and $\{\chi_1^*,\ldots,\chi_3^*\} = \{\chi_{34},\chi_{35},\chi_{36}\}$,
and $\hat\chi = \chi_{44}$.
Decomposing the restrictions of the ordinary characters in $B_1$
into the basic set $\{\chi'_{17},\chi'_{34},\chi'_{44}\}$
we obtain the matrix of Table~\ref{DecMatB1mod2}, where the 
characters belonging to the basic set are indicated by bold face.
This shows that we have 
$$ \begin{array}{rclrclrclrcl}
\chi_1&=&\chi_{37},&\chi_2&=&\chi_{17},&
\chi_3&=&\chi_{49},&\chi_4&=&\chi_{45},\\
\delta_1&=&1,&\delta_2&=&-1,&
\delta_3&=&-1,&\delta_4&=&1,\\
\end{array}
$$
that is, by \cite[Lemma 11.6]{ErdmannSemidihedral} the decomposition matrix of
$B_1$ is as given in Table~\ref{DecMatB1mod2}.
\begin{table}
\[
\begin{array}{r|r|rrr}  \hline
\chi & \chi(1) 
& \Phi_1 & \Phi_2 & \Phi_3 \\ \hline
{\bf 17}&  214\,016 & 1 & . & . \rule[ 0pt]{0pt}{ 13pt} \\
{\bf 34}& 1\,361\,920 & . & 1 & . \\
35& 1\,361\,920 & . & 1 & . \\
36& 1\,361\,920 & . & 1 & . \\
37& 1\,575\,936 & 1 & 1 & . \\
{\bf 44}& 2\,985\,984 & . & . & 1 \\
45& 3\,200\,000 & 1 & . & 1 \\
49& 4\,561\,920 & 1 & 1 & 1   
\rule[- 7pt]{0pt}{ 5pt} \\ \hline
\end{array}
\]
\caption{The decomposition matrix of the 
$2$-block~$B_1$ of $\mathsf{HN}$}
\label{DecMatB1mod2}
\end{table}

\subsection{Proof for the principal block $B_0$}\label{subsec:B0mod2} 
Our proof for the principal block relies on the tried and tested formula
of applying the MeatAxe in conjunction with condensation. Owing to the
modest dimensions of the first few irreducible representations of $G$,
the first five nontrivial ones are available through \cite{web:ATLAS}.
This way we obtain $132_1$, $132_2$, $760$, $2\,650_1$, and $2\,650_2$;
where $2\,650_1$ yields $2\,650_2$ by taking the contragredient. Also
available from the same source is the smallest non-trivial irreducible
representation of the sporadic Baby Monster group $\mathsf{B}$ in
characteristic 2. The restriction of this 4370-dimensional module
yields the composition factors
\[ 4\,370\smash\downarrow_G = 3\,344 + 760 + 132_1 + 132_2,\] 
thus adding $3\,344$
to our collection of simple modules we have available for further
computations. Also owing to the small size of these modules, we may compute
their Brauer characters explicitly by lifting eigenvalues,
where the necessary straight line programmes to restrict from 
$\mathsf{B}$ to $\mathsf{HN}$ and to find conjugacy class representatives 
are also available at \cite{web:ATLAS}.

The eighth Brauer character in the principal block may also be computed
with a plain application of the MeatAxe: Chopping the tensor product
$132_1 \otimes 132_2$ gives the composition factors 
\[ 132_1 \otimes 132_2 = 15\,904 + 2\times 760\]
and lets us derive the Brauer character
of $15\,904$ from the tensor product of the Brauer characters.
As the direct analysis of even larger tensor products with the MeatAxe is
considerably slowed and ultimately rendered infeasible 
due to the increased module dimensions, we 
determine the remaining nine simple modules and their Brauer characters with
the help of condensation.

As our condensation subgroup $K$ we choose a Sylow 3-subgroup of order 243
of the alternating group on 12 letters which in turn is the largest maximal
subgroup of $G$, and condense with the corresponding trace idempotent $e$;
here, a straight line programme to find the largest maximal subgroup of $G$
is available in \cite{web:ATLAS}, and Sylow subgroups of permutation
groups can be computed in \textsf{GAP}.
Unfortunately, for our choice of condensation subgroup, there are at
present no computationally efficient methods available to deal with the
generation problem (see \cite{tackling}). Here we arbitrarily choose
a few elements $g_i\in G$, and have to allow for the possibility that the
corresponding \emph{condensed elements} $eg_ie$ generate 
a proper subalgebra $\Ccal$ of the condensed group algebra $eFGe$.
We apply the condensation of tensor products technique
described in \cite{art:newtensorcondense, LuxWie}.

\begin{table}
 \[
 \begin{array}{r|rrrrrrrrrrrrrrrrr}
  \hline &
  \begin{turn}{90}1\end{turn}&
  \begin{turn}{90}$132_1$\end{turn}&
  \begin{turn}{90}$132_2$\end{turn}&
  \begin{turn}{90}760\end{turn}&
  \begin{turn}{90}$2\,650_1$\end{turn}&
  \begin{turn}{90}$2\,650_2$\end{turn}&
  \begin{turn}{90}3\,344\end{turn}&
  \begin{turn}{90}$132_1 \otimes 132_2$\end{turn}&
  \begin{turn}{90}$132_1 \otimes 760$\end{turn}&
  \begin{turn}{90}$132_2\otimes 760$\end{turn}&
  \begin{turn}{90}$760\otimes 760$\end{turn}&
  \begin{turn}{90}$132_1 \otimes 2\,650_1$\end{turn}&
  \begin{turn}{90}$132_1 \otimes 3\,344$\end{turn}&
  \begin{turn}{90}$132_2 \otimes 3\,344$\end{turn}&
  \begin{turn}{90}$760 \otimes 2\,650_1$\end{turn}&
  \begin{turn}{90}$760 \otimes 3\,344$\end{turn}&
  \begin{turn}{90}$2\,650_1 \otimes 2\,650_2$\end{turn}\\\hline
  k1& 1& .& .& .& .& .& .& .& 2& 2&  20& 6& 6& 6&  40&  12&  46 \\
 k2_1 & .& 1& .& .& .& .& .& .& 2& 2& 6& 3& 4& 4&  22& 6&  23 \\
 k2_2 & .& .& 1& .& .& .& .& .& 2& 2& 6& 4& 4& 4&  22& 6&  23 \\
 k6_1 & .& .& .& .& 1& .& .& .& .& .& 8& .& .& .& 5& 2& 8 \\
 k6_2 & .& .& .& .& .& 1& .& .& .& .& 8& 1& .& .& 5& 2& 8 \\
 k12 & .& .& .& 1& .& .& .& 1& .& .& 3& 1& 2& 2&  10& 5&  11 \\
 k28 & .& .& .& .& .& .& 1& .& .& .&  10& .& .& .& 2& 4& 8 \\
 k62 & .& .& .& .& .& .& .& 1& .& .& 4& 1& 2& 2&  12& 4&  10 \\
 k118_1 & .& .& .& .& .& .& .& .& .& 1& .& 1& .& .& .& .& 2 \\
 k118_2 & .& .& .& .& .& .& .& .& 1& .& .& .& .& .& .& .& 2 \\
 k160_1 & .& .& .& .& .& .& .& .& .& .& .& 1& 1& 1& 7& 1& 5 \\
 k160_2 & .& .& .& .& .& .& .& .& .& .& .& 1& 1& 1& 6& 1& 5 \\
 k164_1 & .& .& .& .& .& .& .& .& 2& .& .& .& 1& 2& 1& 2& 4 \\
 k164_2 & .& .& .& .& .& .& .& .& .& 2& .& 1& 2& 1& 1& 2& 4 \\
 k706 & .& .& .& .& .& .& .& .& .& .& .& 1& .& .& 4& .& 3 \\
 k922 & .& .& .& .& .& .& .& .& .& .& 2& .& 1& 1& 2& 2& 2 \\
 k6\,344 & .& .& .& .& .& .& .& .& .& .& .& .& .& .& .& .& 1 \\\hline
 k908 & .&.&.&.&.&.&.& .&  .&  .&  .&  .&  .&  .&  .&  2&  . \\
 k5\,652 & .&.&.&.&.&.&.& .&  .&  .&  .&  .&  .&  .&  .&  1&  . \\
 k14\,072 & .&.&.&.&.&.&.& .&  .&  .&  .&  .&  .&  .&  .&  .&  1 \\\hline
 \end{array}
 \]
 \caption{Condensation results for $p=2$ for $\mathsf{HN}$}
 \label{tab:multmod2}
\end{table}

The results are presented in Table~\ref{tab:multmod2}, in which we
display the multiplicities of the composition factors of the condensed
tensor products given, where we denote the simple $\Ccal$-modules by
their dimension, preceded by the letter `$k$'. We also condense the
known simple $FG$-modules of the principal block individually in order to
match them to simple $\Ccal$-modules. Note that we condense the tensor
product $132_1 \otimes 132_2$, even though we were able to compute its
decomposition directly: it is more efficient to condense the tensor
product than to deal with the module $15\,904$ alone.

Our next task is to determine which simple $\Ccal$-modules are
composition factors of simple $eFGe$-modules coming from simple
$FG$-modules of its principal block. To this end, it turns out that
the restrictions of the ordinary characters indicated by bold face in
Table~\ref{tab:decB0mod2} are a basic set of Brauer characters for the
principal block. A basic set of Brauer characters for block $B_1$ was
already given in Table~\ref{DecMatB1mod2}, and adding $\chi_{46}'$ we
obtain a complete basic set of Brauer characters.

Then, by decomposing the Brauer characters of the tensor
products considered into this basic set, we see that only the last two
contain composition factors which lie outside of the principal block.
As a manner of speaking, we say that a simple $\Ccal$-module belongs
to a block of $FG$ if it is a composition factor of the restriction to
$\Ccal$ of a condensed simple $FG$-module lying in that block.

By computing the multiplicity of the trivial character of $K$ in the
restrictions to $K$ of the basic set characters and the known Brauer
characters, we determine the dimensions of the corresponding condensed
modules. Of course, these are upper bounds for the dimensions of their
simple $\Ccal$-constituents.

Also, decomposing the characters of the last two tensor products in
the basic set of $B_0$ and the irreducible Brauer characters of the other
blocks, gives that
$214\,016$ occurs twice, and $1\,361\,920$ occurs once as a constituent of 
$760\otimes 3\,344$. The single defect zero module $3\,424\,256$ 
appears once in $2\,650_1 \otimes 2\,650_2$.

We compute the dimensions of the condensed modules whose
Brauer characters constitute the basic set. From this it is immediate that
$k14\,072$ is a composition factor of the condensed defect zero module
$3\,424\,256$ restricted to $\Ccal$: it is simply to large to fit into any
of the other modules. As the dimension of the $\Ccal$-module is equal to
the dimension computed of $3\,424\,256$ condensed, we conclude that
$k14\,072$ is the irreducible restriction of the latter module, and
therefore belongs to the defect zero block.

The simple $1\,361\,920$-dimensional modules of $B_1$ condense to modules of
dimension $5\,652$. By summing up the dimensions of the composition factors
by their multiplicities as given in the second to last column of
Table~\ref{tab:multmod2}, it is plain to see 
that $k5\,652$ appears in one of the condensed simple modules of the same
dimension, as the sum obtained is less.
Thus $k5\,652$ too is an irreducible restriction of a condensed simple
module lying in $B_1$ and therefore belongs to $B_1$.

For $k908$ we verify computationally that it also originates from a module
outside of the principal block: employing the MeatAxe and peakword
techniques, we find a uniserial
submodule of $760 \otimes 3\,344$ condensed
with head and socle isomorphic to $k908$ and containing $k5\,652$
as the only additional composition factor. Therefore there is a
non-trivial extension
of $k5\,652$ by $k908$, and as $k5\,652$ belongs to $B_1$ so must $k908$.

Thus all appearing $\Ccal$-modules, except $k908$, $k5\,652$, and
$k14\,072$, belong to the principal block. Hence the matrix $A$
comprising the first $17$ rows of Table~\ref{tab:multmod2} describes the
condensed modules of the $B_0$-components of the modules given at its
top. It is immediate that $A$ has full rank. Therefore, as $B_0$ has
$l=17$ modular characters, by \cite[Proposition~1]{HenkeHissMueller}
condensation with respect to our chosen subgroup $K$ defines a
faithful functor from $\mathsf{mod}$-$B_0$ to $\mathsf{mod}$-$eB_0e$.
Furthermore, by \cite[Theorem~4.1]{LuxNoeskeRyba} taking as $b$ the
$B_0$-components of the Brauer characters afforded by the modules indexing
the columns of Table~\ref{tab:multmod2}, 
the product $bA^{-1}$ gives a basic set of Brauer atoms for $B_0$.

As the possibility that $\Ccal$ is a proper subalgebra of $eFGe$
remains, the degrees of the Brauer atoms are lower bounds for the
degrees of the irreducible 
Brauer characters. Therefore in order to verify that the
Brauer atoms are in fact the Brauer characters, we need to show that
the simple $FG$-modules possess these dimensions. This is achieved by
constructing bases for almost all simple modules through uncondensation:
With the help of the MeatAxe we can determine the socle series of the
condensed tensor products considered, and may therefore seek a suitably
small submodule $U$ of the condensed module $Ve$ which contains a
condensed version of the simple module $S$ we would like to construct.
The subspace embedding of $Ve$ into $V$ therefore embeds a basis of
$U$ into $V$, and with standard MeatAxe functionality we may construct
a basis of the smallest $G$-invariant subspace of $V$ containing the
embedded basis of $U$. Ideally, we have $U =Se$, and we are done after
having constructed a basis of $S$. However, if $Se$ does not lie in the
socle of $V$, we need to iterate this procedure by working along an
ascending composition series of $U$.

Recall that we have already chopped up the tensor product $132_1 \otimes
132_2$ completely, hence we are done with that. Within the condensed
tensor products $132_1 \otimes 760$, $760 \otimes 760$, $132_1\otimes
2\,650_1$, and $132_1\otimes 3\,344$, respectively, we find suitable
submodules whose socle series we give in Table~\ref{tab:socseries}.
\begin{table}
 \centering
 \begin{tabular}{cccc}
  \boxed{
  \begin{tabular}{c}
   $k118_2$ \\ $k2_2$ \\ $k1$ \\ $k2_1 \oplus k164_1$
  \end{tabular}
  }
  &
  \boxed{
  \begin{tabular}[h]{c}
   $k922$\\ $k1 \oplus k62$ \\ $k2_1 \oplus k2_2 \oplus k28$ \\
   $k1 \oplus k1$ \\ $k6_1 \oplus k6_2$ \\ $k6_1 \oplus k6_2$ \\
   $k28$
  \end{tabular}
  }
  &
  \boxed{
  \begin{tabular}[h]{c}
   $k706$ \\ $k1$ \\ $k2_2$
  \end{tabular}
  }
  &
  \boxed{
  \begin{tabular}[h]{c}
   $k2_1 \oplus k12 \oplus k160_2$\\
   $k164_2$
  \end{tabular}
  }
 \end{tabular}
 \caption{Socle series of select submodules of some tensor products}
 \label{tab:socseries}
\end{table}
Note that for every pair of modules conjugate under the outer
automorphism, which of course is seen on the level of characters,
it is sufficient to construct a basis of only one member. Therefore
the four modules of Table~\ref{tab:socseries} allow us to verify the
dimensions of almost all Brauer characters except the one corresponding
to $k6\,344$. For the practical realisation of these computations, it is
important to note that we do not need to construct a representation of
the action of $G$ on any of the huge tensor product spaces: using the
natural isomorphism $M\otimes N \cong \Hom_F(M^*,N)$ (see \cite[Lemma
2.6]{LuxWiegelmannSoc}) we only need to work with the available
representations of the tensor factors.

With only the last Brauer atom of degree $1\,556\,136$ left to prove
to be a Brauer character, we can avoid the unwieldy construction of a
basis for this module: We additionally condense the tensor product $V :=
2\,650_1 \otimes 3\,344$ instead. Let $S$ denote the simple $FG$-module
such that its Brauer character $\varphi$ contains the Brauer atom
$\alpha$ of degree $1\,556\,136$ as a summand. Since $\varphi$ is a
constituent of the basic set character $\chi'_{39}$, and the latter has
degree $2\,031\,480$, we conclude that $\varphi$ contains $\alpha$ as a
summand with multiplicity precisely one.

A calculation with the MeatAxe gives that head and socle of
$Ve\smash\downarrow_\Ccal$ are isomorphic, and have structure $$ k160_1
\oplus k6\,344 \oplus k12\,288 .$$ By the above we have $43\,416_1 e =
k160_1$ and $2\,985\,984 e = k12\,288$, and $Se\smash\downarrow_\Ccal$
has $k6\,344$ as a composition factor. Recall that $e$ is a faithful
condensation idempotent, and that the Brauer atom $\alpha$ alone already
accounts for the condensed dimension $6\,344$, that is $\varphi=\alpha$
is equivalent to $Se=k6\,344$, but due to the generation problem this
does not immediately follow from the above statement, and neither it
is clear how socle and head of $Ve$ and thus of $V$ look like. But,
with the possible exception of $\varphi$, we have proved that all
simple $eFGe$-modules restrict irreducibly to $\Ccal$, and that two such
$eFGe$-modules are isomorphic if their restrictions to $\Ccal$ are. Hence
we conclude that both the head and the socle of $V$ have at most the
constituents $43\,416_1$, $2\,985\,984$ and $S$; in particular, the
constituents $132_{1/2}$ and $760$ do not occur.

We now additionally choose another condensation subgroup, namely
the normal subgroup $5_+^{1+4}$ of the sixth maximal subgroup
$M_6:=5_+^{1+4}\colon 2_-^{1+4}.5.4$. To distinguish this condensation
from the previous one, we denote the trace idempotent of $5_+^{1+4}$
by $e'$ and prefix the simple $e'FGe'$-modules by the letter '$c$'.
By \cite[Theorem~2.7]{tackling} we obtain a generating system of $e'FGe'$ of
$128$ elements, if we condense the $126$ non-identity double coset
representatives of $M_6\backslash G/M_6$ and the two generators for $M_6$,
obtained through a straight line programme available at \cite{web:ATLAS}.

As we know all irreducible Brauer characters except $\varphi$, we can
determine the dimensions of their condensed modules with respect to
$e'$. This shows that $e'$ is not faithful; indeed, precisely the simple
modules $132_{1/2}$ and $760$ condense to the zero module. We
identify the possible composition factors of the $2\,772$-dimensional
condensed module $Ve'$ by applying the MeatAxe. We again find that head
and socle of $Ve'$ are isomorphic, and have structure $$ c8 \oplus
c1\,152 \oplus c960 ,$$ where $43\,416_1 e' = c8$, $2\,985\,984 e' =
c960$ and $Se'= c1\,152$. Again, a consideration of condensed dimensions
shows that $\alpha$ alone accounts for the dimension $1\,152$, which
implies that $\varphi-\alpha$ is a sum of multiples of $132_{1/2}$ and
$760$.

Moreover, as we look at $Ve'$ as an $e'FGe'$-module, uncondensing (in
theory, not in practice!) the submodule of $Ve'$ isomorphic to $c1\,152$
yields a submodule $U$ of $V$ having a simple head isomorphic to $S$,
whose radical has only constituents condensing to the zero module, that
is amongst $132_{1/2}$ and $760$. But since none of the latter occur in
the socle of $V$ we conclude that $U$ is isomorphic to $S$, in other
words $S$ occurs in the socle of $V$, and similarly in the head of $V$
as well.

We now return to our established faithful condensation setup and
reconsider $Ve$: Assume that $\varphi\neq\alpha$. Then, as $\alpha$ is a
summand of $\varphi$ with multiplicity one, hence $k6\,344$ occurs only
once as a constituent of $Se\smash\downarrow_\Ccal$, we conclude that the
head or the socle of $Ve\smash\downarrow_\Ccal$ (or both) contain a
condensed module of at least one of $132_{1/2}$ and $760$, which is not
the case, and thus is a contradiction. Hence we conclude that the last
Brauer atom $\alpha$ coincides with $\varphi$, thus is an irreducible
Brauer character.

The decomposition matrix of the principal block is given in
Table~\ref{tab:decB0mod2}.

\begin{table}
{\tiny
\hspace*{-7em} $
 \begin{array}{r|r|rrrrrrrrrrrrrrrrr}
  \hline
\chi & \chi(1) & \varphi_1 & \varphi_2 & \varphi_3 & \varphi_4 & \varphi_5 & \varphi_6 &
\varphi_7 & \varphi_8 & \varphi_9 & \varphi_{10} & \varphi_{11} & \varphi_{12} &
\varphi_{13} & \varphi_{14} & \varphi_{15} & \varphi_{16} & \varphi_{17}\\
\hline
{\bf 1}&       1  & 1& .& .& .& .& .& .& .& .& .& .& .& .& .& .& .& . \\
{\bf 2}&      133&  1& .& 1& .& .& .& .& .& .& .& .& .& .& .& .& .& . \\
{\bf 3}&      133&  1& 1& .& .& .& .& .& .& .& .& .& .& .& .& .& .& . \\
{\bf 4}&      760&  .& .& .& .& .& 1& .& .& .& .& .& .& .& .& .& .& . \\
{\bf 5}&   3\,344&  .& .& .& .& .& .& 1& .& .& .& .& .& .& .& .& .& . \\
{\bf 6}&   8\,778&  2& .& 1& 1& 1& .& 1& .& .& .& .& .& .& .& .& .& . \\
7&   8\,778&  2& 1& .& 1& 1& .& 1& .& .& .& .& .& .& .& .& .& . \\
8&   8\,910&  2& 1& 1& 1& 1& .& 1& .& .& .& .& .& .& .& .& .& . \\
9&   9\,405&  1& .& .& 1& 1& 1& 1& .& .& .& .& .& .& .& .& .& . \\
{\bf 10}&  16\,929&  1& 1& 1& .& .& 1& .& 1& .& .& .& .& .& .& .& .& . \\
{\bf 11}&  35\,112&  .& .& .& .& .& 1& .& .& .& .& .& .& .& 1& .& .& . \\
{\bf 12}&  35\,112&  .& .& .& .& .& 1& .& .& .& .& .& .& 1& .& .& .& . \\
{\bf 13}&  65\,835&  1& 2& 1& .& .& .& .& .& 1& .& .& .& .& 1& .& .& . \\
{\bf 14}&  65\,835&  1& 1& 2& .& .& .& .& .& .& 1& .& .& 1& .& .& .& . \\
{\bf 15}&  69\,255&  3& 2& 2& 1& 1& 1& 1& 1& .& .& 1& .& .& .& .& .& . \\
{\bf 16}&  69\,255&  3& 2& 2& 1& 1& 1& 1& 1& .& .& .& 1& .& .& .& .& . \\
{\bf 18}& 267\,520&  6& 2& 2& 2& 2& 1& 2& 2& .& .& .& .& .& .& .& 1& . \\
19& 270\,864&  6& 2& 2& 2& 2& 1& 3& 2& .& .& .& .& .& .& .& 1& . \\
{\bf 20}& 365\,750&  8& 3& 3& 2& 2& 3& 1& 1& .& .& 1& 1& 1& 1& 1& .& . \\
21& 374\,528&  6& 2& 2& 1& 1& 2& .& .& .& .& 1& 1& 2& 1& 1& .& . \\
22& 374\,528&  6& 2& 2& 1& 1& 2& .& .& .& .& 1& 1& 1& 2& 1& .& . \\
23& 406\,296&  6& 4& 4& .& .& 1& .& 2& .& .& 1& 1& 1& 1& .& 1& . \\
24& 653\,125&  13& 7& 7& 2& 2& 3& 2& 4& .& .& 2& 2& .& .& 1& 1& . \\
{\bf 25}& 656\,250&  12& 6& 6& 1& 2& 3& 1& 3& .& .& 1& 2& 1& 1& 1& 1& . \\
26& 656\,250&  12& 6& 6& 2& 1& 3& 1& 3& .& .& 2& 1& 1& 1& 1& 1& . \\
27& 718\,200&  14& 8& 9& 2& 2& 2& 2& 4& .& 1& 2& 2& 1& .& 1& 1& . \\
28& 718\,200&  14& 9& 8& 2& 2& 2& 2& 4& 1& .& 2& 2& .& 1& 1& 1& . \\
29&1\,053\,360&  26& 12& 12& 6& 6& 7& 5& 6& .& .& 3& 3& 1& 1& 2& 1& . \\
30&1\,066\,527&  19& 11& 11& 2& 2& 5& 1& 5& 1& .& 3& 3& 1& 2& 2& 1& . \\
31&1\,066\,527&  19& 11& 11& 2& 2& 5& 1& 5& .& 1& 3& 3& 2& 1& 2& 1& . \\
32&1\,185\,030&  20& 12& 12& 2& 2& 5& 1& 4& 1& 1& 3& 3& 3& 3& 2& 1& . \\
33&1\,354\,320&  24& 14& 14& 2& 2& 6& 1& 8& .& .& 4& 4& 1& 1& 2& 2& . \\
38&1\,625\,184&  30& 16& 16& 4& 4& 7& 4& 10& .& .& 4& 4& 1& 1& 2& 3& . \\
{\bf 39}&2\,031\,480&  14& 8& 8& 2& 2& 2& 2& 4& .& .& 2& 2& .& .& .& 1& 1 \\
40&2\,375\,000&  20& 8& 8& 4& 4& 6& 2& 2& .& .& 3& 3& 2& 2& 2& .& 1 \\
41&2\,407\,680&  50& 26& 26& 8& 8& 13& 6& 14& .& .& 7& 7& 2& 2& 4& 3& . \\
42&2\,661\,120&  48& 28& 28& 5& 5& 8& 4& 14& 1& 1& 7& 7& 2& 2& 4& 4& . \\
43&2\,784\,375&  25& 14& 14& 3& 3& 5& 2& 5& 1& 1& 4& 4& 2& 2& 2& 1& 1 \\
47&3\,878\,280&  48& 27& 27& 7& 7& 10& 6& 12& 1& 1& 7& 7& 3& 3& 3& 3& 1 \\
48&4\,156\,250&  52& 29& 29& 7& 7& 11& 5& 13& 1& 1& 8& 8& 3& 3& 4& 3& 1 \\
50&4\,809\,375&  65& 36& 36& 9& 9& 14& 7& 17& 1& 1& 10& 10& 3& 3& 5& 4& 1 \\
51&5\,103\,000&  68& 38& 38& 8& 8& 14& 6& 18& 1& 1& 10& 10& 4& 4& 5& 5& 1 \\
52&5\,103\,000&  68& 38& 38& 8& 8& 14& 6& 18& 1& 1& 10& 10& 4& 4& 5& 5& 1 \\ 
53&5\,332\,635&  71& 41& 41& 8& 8& 16& 5& 17& 2& 2& 11& 11& 6& 6& 6& 4& 1 \\ 
54&5\,878\,125&  81& 47& 47& 9& 9& 16& 7& 22& 2& 2& 12& 12& 5& 5& 6& 6& 1 \\
\hline
 \end{array}
 $}
\hspace*{0.5em}
 \caption{The decomposition matrix of the principal 
2-block $B_0$ of $\mathsf{HN}$}
 \label{tab:decB0mod2}
\end{table}

\subsection{The Automorphism Group in Characteristic 2}
The $2$-modular characters of $G.2$ are distributed over three blocks,
each covering precisely one of the $2$-blocks of $G$: the principal
block $B_0$ of defect $d=15$, having $k=63$ ordinary and $l=12$ modular
characters, the block $B_1$ of defect $d=5$, having $k=13$ ordinary and
$l=3$ modular characters, and the block $B_2$ of defect $d=1$, having
$k=2$ ordinary and $l=1$ modular characters.

The task of determining the decomposition numbers of $G.2$ in
characteristic 2 is easy: every Brauer character of $G$ invariant under
the outer automorphism gives just one Brauer character of $G.2$, while
the pairs of Brauer characters of $G$ whose members are interchanged by
the outer automorphism induce to one Brauer character of $G.2$. Therefore,
using the results of the previous section we can immediately give the
three decomposition matrices in Tables~\ref{tab:decmatB0HN2mod2} and
\ref{tab:decmatB12HN2mod2}.

\begin{table}
\tiny
 \[
 {
 \begin{array}{r|r|rrrrrrrrrrrr} \hline
\chi &\chi(1) & \Phi_1 & \Phi_2 & \Phi_3 & \Phi_4 & \Phi_5 & \Phi_6 & \Phi_7
  & \Phi_8 & \Phi_9 & \Phi_{10} & \Phi_{11} & \Phi_{12} \\ \hline 
  1 &      1 & 1& .& .& .& .& .& .& .& .& .& .& . \\
  2 &      1 & 1& .& .& .& .& .& .& .& .& .& .& . \\
  3 &    266 & 2& 1& .& .& .& .& .& .& .& .& .& . \\
  4 &    760 & .& .& .& 1& .& .& .& .& .& .& .& . \\
  5 &    760 & .& .& .& 1& .& .& .& .& .& .& .& . \\
  6 &   3\,344 & .& .& .& .& 1& .& .& .& .& .& .& . \\
  7 &   3\,344 & .& .& .& .& 1& .& .& .& .& .& .& . \\
  8 &  17\,556 & 4& 1& 2& .& 2& .& .& .& .& .& .& . \\
  9 &   8\,910 & 2& 1& 1& .& 1& .& .& .& .& .& .& . \\
 10 &   8\,910 & 2& 1& 1& .& 1& .& .& .& .& .& .& . \\
 11 &   9\,405 & 1& .& 1& 1& 1& .& .& .& .& .& .& . \\
 12 &   9\,405 & 1& .& 1& 1& 1& .& .& .& .& .& .& . \\
 13 &  16\,929 & 1& 1& .& 1& .& 1& .& .& .& .& .& . \\
 14 &  16\,929 & 1& 1& .& 1& .& 1& .& .& .& .& .& . \\
 15 &  70\,224 & .& .& .& 2& .& .& .& .& 1& .& .& . \\
 16 & 131\,670 & 2& 3& .& .& .& .& 1& .& 1& .& .& . \\
 17 & 138\,510 & 6& 4& 2& 2& 2& 2& .& 1& .& .& .& . \\
 20 & 267\,520 & 6& 2& 2& 1& 2& 2& .& .& .& .& 1& . \\
 21 & 267\,520 & 6& 2& 2& 1& 2& 2& .& .& .& .& 1& . \\
 22 & 270\,864 & 6& 2& 2& 1& 3& 2& .& .& .& .& 1& . \\
 23 & 270\,864 & 6& 2& 2& 1& 3& 2& .& .& .& .& 1& . \\
 24 & 365\,750 & 8& 3& 2& 3& 1& 1& .& 1& 1& 1& .& . \\
 25 & 365\,750 & 8& 3& 2& 3& 1& 1& .& 1& 1& 1& .& . \\
 26 & 749\,056 & 12& 4& 2& 4& .& .& .& 2& 3& 2& .& . \\
 27 & 406\,296 & 6& 4& .& 1& .& 2& .& 1& 1& .& 1& . \\
 28 & 406\,296 & 6& 4& .& 1& .& 2& .& 1& 1& .& 1& . \\
 29 & 653\,125 & 13& 7& 2& 3& 2& 4& .& 2& .& 1& 1& . \\
 30 & 653\,125 & 13& 7& 2& 3& 2& 4& .& 2& .& 1& 1& . \\
 31 &1\,312\,500 & 24& 12& 3& 6& 2& 6& .& 3& 2& 2& 2& . \\
 32 &1\,436\,400 & 28& 17& 4& 4& 4& 8& 1& 4& 1& 2& 2& . \\
 33 &1\,053\,360 & 26& 12& 6& 7& 5& 6& .& 3& 1& 2& 1& . \\
 34 &1\,053\,360 & 26& 12& 6& 7& 5& 6& .& 3& 1& 2& 1& . \\
 35 &2\,133\,054 & 38& 22& 4& 10& 2& 10& 1& 6& 3& 4& 2& . \\
 36 &1\,185\,030 & 20& 12& 2& 5& 1& 4& 1& 3& 3& 2& 1& . \\
 37 &1\,185\,030 & 20& 12& 2& 5& 1& 4& 1& 3& 3& 2& 1& . \\
 38 &1\,354\,320 & 24& 14& 2& 6& 1& 8& .& 4& 1& 2& 2& . \\
 39 &1\,354\,320 & 24& 14& 2& 6& 1& 8& .& 4& 1& 2& 2& . \\
 45 &1\,625\,184 & 30& 16& 4& 7& 4& 10& .& 4& 1& 2& 3& . \\
 46 &1\,625\,184 & 30& 16& 4& 7& 4& 10& .& 4& 1& 2& 3& . \\
 47 &2\,031\,480 & 14& 8& 2& 2& 2& 4& .& 2& .& .& 1& 1 \\
 48 &2\,031\,480 & 14& 8& 2& 2& 2& 4& .& 2& .& .& 1& 1 \\
 49 &2\,375\,000 & 20& 8& 4& 6& 2& 2& .& 3& 2& 2& .& 1 \\
 50 &2\,375\,000 & 20& 8& 4& 6& 2& 2& .& 3& 2& 2& .& 1 \\
 51 &2\,407\,680 & 50& 26& 8& 13& 6& 14& .& 7& 2& 4& 3& . \\
 52 &2\,407\,680 & 50& 26& 8& 13& 6& 14& .& 7& 2& 4& 3& . \\
 53 &2\,661\,120 & 48& 28& 5& 8& 4& 14& 1& 7& 2& 4& 4& . \\
 54 &2\,661\,120 & 48& 28& 5& 8& 4& 14& 1& 7& 2& 4& 4& . \\
 55 &2\,784\,375 & 25& 14& 3& 5& 2& 5& 1& 4& 2& 2& 1& 1 \\
 56 &2\,784\,375 & 25& 14& 3& 5& 2& 5& 1& 4& 2& 2& 1& 1 \\
 63 &3\,878\,280 & 48& 27& 7& 10& 6& 12& 1& 7& 3& 3& 3& 1 \\
 64 &3\,878\,280 & 48& 27& 7& 10& 6& 12& 1& 7& 3& 3& 3& 1 \\
 65 &4\,156\,250 & 52& 29& 7& 11& 5& 13& 1& 8& 3& 4& 3& 1 \\
 66 &4\,156\,250 & 52& 29& 7& 11& 5& 13& 1& 8& 3& 4& 3& 1 \\
 69 &4\,809\,375 & 65& 36& 9& 14& 7& 17& 1& 10& 3& 5& 4& 1 \\
 70 &4\,809\,375 & 65& 36& 9& 14& 7& 17& 1& 10& 3& 5& 4& 1 \\
 71 &5\,103\,000 & 68& 38& 8& 14& 6& 18& 1& 10& 4& 5& 5& 1 \\
 72 &5\,103\,000 & 68& 38& 8& 14& 6& 18& 1& 10& 4& 5& 5& 1 \\
 73 &5\,103\,000 & 68& 38& 8& 14& 6& 18& 1& 10& 4& 5& 5& 1 \\
 74 &5\,103\,000 & 68& 38& 8& 14& 6& 18& 1& 10& 4& 5& 5& 1 \\
 75 &5\,332\,635 & 71& 41& 8& 16& 5& 17& 2& 11& 6& 6& 4& 1 \\
 76 &5\,332\,635 & 71& 41& 8& 16& 5& 17& 2& 11& 6& 6& 4& 1 \\
 77 &5\,878\,125 & 81& 47& 9& 16& 7& 22& 2& 12& 5& 6& 6& 1 \\
 78 &5\,878\,125 & 81& 47& 9& 16& 7& 22& 2& 12& 5& 6& 6& 1 \\ \hline  
 \end{array}}
 \]
 \caption{The decomposition matrix of the principal 2-block $B_0$ 
          of $\mathsf{HN}.2$}
 \label{tab:decmatB0HN2mod2}
\end{table}

\begin{table}
 \[
 \begin{array}{r|r|rrr|r} \hline
\chi &\chi(1) & \Phi_1 & \Phi_2 & \Phi_3 & \Phi_1 \\ \hline 
18&  214\,016 &  1&  .&  . &. \\
19&  214\,016 &  1&  .&  . &. \\
40& 1\,361\,920 & .&  1&  . &. \\
41& 1\,361\,920 &  .&  1&  . &. \\
42& 2\,723\,840 &  .&  2&  . &. \\
43& 1\,575\,936 &  1&  1&  . &. \\
44& 1\,575\,936 &  1&  1&  . &. \\
57& 2\,985\,984 &  .&  .&  1 &. \\
58& 2\,985\,984 &  .&  .&  1 &. \\
59& 3\,200\,000 &  1&  .&  1 &. \\
60& 3\,200\,000 &  1&  .&  1 &. \\
67& 4\,561\,920 &  1&  1&  1 &. \\
68& 4\,561\,920 &  1&  1&  1 &. \\ \hline
61&  3\,424\,256 & .&.&.& 1 \\
62&  3\,424\,256 & .&.&.& 1 \\ \hline
 \end{array}
 \]
 \caption{The decomposition matrices of the 2-blocks $B_1$ and $B_2$ of 
          $\mathsf{HN}.2$}
 \label{tab:decmatB12HN2mod2}
\end{table}

\section{Characteristic~$3$}
\label{sec:mod3}
In characteristic 3, the simple group $G$ possesses three
blocks of positive defect: the principal block $B_0$ of defect
$d=6$, having $k=33$ ordinary and $l=20$ modular characters,
the block $B_1$ of defect $d=2$, having $k=9$ ordinary and
$l=7$ modular characters, and the block $B_2$ of defect $d=1$,
having $k=3$ ordinary and $l=2$ modular characters. Moreover,
there are nine blocks of defect $0$, consisting of the ordinary
characters $\{\chi_{15},\chi_{16},\chi_{30},\chi_{31},\chi_{44},
\chi_{47},\chi_{51},\chi_{52},\chi_{53}\}$.

The decomposition matrix of the $3$-block $B_2$, reproduced in
Table~\ref{DecMatB2}, follows immediately from the Brauer-Dade theory of
blocks with cyclic defect groups, and has previously been published in
\cite{BrauerBaumBuch}.
\begin{table}
\[
\begin{array}{r|r|rr} \hline
\chi & \chi(1) & \Phi_1 & \Phi_2 \\ \hline
23&  406\,296 &  1 & . \rule[ 0pt]{0pt}{ 13pt} \\
38& 1\,625\,184 &  . & 1 \\
39& 2\,031\,480 &  1 & 1
\rule[- 7pt]{0pt}{ 5pt} \\ \hline
\end{array}
\]
\caption{The decomposition matrix of the
$3$-block $B_2$ of $\mathsf{HN}$}
\label{DecMatB2}
\end{table}

\subsection{Proof for the block~$B_1$, part 1}
\label{subsec:B1mod3part1}
\begin{table}
$$
\begin{array}{r|r|rrrrrrr|r} \hline
\chi & \chi(1) & \Psi_1 & \Psi_2 & \Psi_3 & \Psi_4 & \Psi_5 & \Psi_6 & 
\Psi_7 & \Psi \\ \hline
{\bf 8}&   8\,910 &  1 &   . &   . &   . &   . &   . &   . &  
     3 \rule[ .pt]{.pt}{ 13pt} \\
{\bf 10}&  16\,929 &  . &   1 &   . &   . &   . &   . &   . & 1 \\
{\bf 19}& 270\,864 &  . &   . &   1 &   . &   . &   . &   . & . \\
{\bf 32}&1\,185\,030 &  1 &   1 &   . &   3 &   1 &   1 &   . & 5 \\
{\bf 33}&1\,354\,320 &  2 &   2 &   . &   2 &   2 &   1 &   . & 4 \\
{\bf 37}&1\,575\,936 &  1 &   2 &   1 &   2 &   1 &   1 &   . & 1 \\
43&2\,784\,375 &  2 &   2 &   1 &   5 &   3 &   2 &   . & 5 \\
{\bf 49}&4\,561\,920 &  1 &   2 &   2 &   9 &   6 &   4 &   1 & 3 \\
50&4\,809\,375 &  1 &   3 &   3 &   9 &   5 &   4 &   1 & 4  
\rule[- 7pt]{0pt}{ 5pt} \\ \hline
\end{array} $$
$$ \Psi_1=\Ind_{M_6}^G( \varphi_{2} )\cdot 1_{B_1}, \quad
   \Psi_3=\Ind_{M_6}^G( \varphi_{3} )\cdot 1_{B_1}, \quad
   \Psi_5=(\Ind_{M_6}^G( \varphi_{7} )\cdot 1_{B_1})/2 \rule{0em}{1em} $$
$$ \Psi_2=(\chi_{2} \cdot \chi_{31}) \cdot 1_{B_1}, \quad\!
   \Psi_4=(\chi_{2} \cdot \chi_{47}) \cdot 1_{B_1}, \quad\!
   \Psi_6=(\chi_{5} \cdot \chi_{15}) \cdot 1_{B_1}, \quad\!
   \Psi_7=(\chi_{4} \cdot \chi_{15}) \cdot 1_{B_1} \rule{0em}{1em} $$
\vspace*{0em}
\caption{\label{PS1} A first basic set of projective characters for~$B_1$}
\end{table}

We begin with a basic set of projective characters displayed in
Table~\ref{PS1}, where we also give a basic set of Brauer characters,
comprising of the restrictions of the ordinary characters indicated by
bold face. The origin of the projective characters is given there as
well, where $1_{B_1}$ denotes the associated block idempotent, $M_6 =
5_+^{1+4}\colon2_-^{1+4}.5.4$ is the sixth maximal subgroup of $G$ (see
\cite{ATLAS}), which is a $3'$-subgroup, the character $\varphi_2$ is
the linear character of~$M_6$ of order~$2$, the character $\varphi_3$
is one of the two linear characters of order~$4$, and $\varphi_7$ is
one of two the rational valued characters of degree~$5$. Note that
the character of the $B_1$-component of $\Ind_{M_6}^G( \varphi_{7} )$
decomposes into the ordinary characters of $B_1$ with even coefficients,
hence dividing by $2$ still yields a projective character. The remaining
projective characters are products of irreducible characters of~$G$, one
factor of which is a character of defect zero.

\medskip
We now pursue a strategy first employed in \cite{MueHabil}. Let
$M'_6$ denote the unique subgroup of $M_6$ of index~$2$. Hence the
$RG$-permutation module on the cosets of~$M'_6$ is a projective
$RG$-lattice. Put \[E := \mbox{\rm End}_{RG}(\Ind_{M'_6}^G(R)\cdot
1_{B_1}) ;\] note that $(Q,R,F)$ is a splitting system for $E$.
The $B_1$-component of the associated permutation character
$\Psi:=\Ind_{M_6}^G(\varphi_1 + \varphi_2)$, where $\varphi_1$ is the
trivial character, decomposes into the ordinary characters of $B_1$
as given in Table~\ref{PS1}. In particular, $Q \otimes_R E$ has eight
irreducible characters.

We consider~$\bar{E} := F \otimes_R E$, and proceed to determine its
regular representation. In practice, a straight line programme to find
$M_6$ is available in \cite{web:ATLAS}. Using the facilities provided
by \textsf{GAP} it is then easy to find $M'_6$, and the permutation
representation of $G$ on the 273\,030\,912 cosets in $M'_6\backslash G$,
realised as a $G$-orbit of vectors in a suitable matrix representation
of $G$. Then the regular representation of~$\bar{E}$ is obtained by a
direct condense technique, using the \textsf{orb} package \cite{orb};
for details see \cite{MueHabil}.

It turns out that $\bar{E}$ has five simple modules $3a,1a,1b,2c,2a$ of
the respective dimensions, and the Cartan matrix~$C$ of~$\bar{E}$ is
found to be as given in Table~\ref{EPims}. Let $D \in \mathbb{N}_0^{8
\times 5}$ be the decomposition matrix of~$E$. It is easy to check that
the matrix equation $D^tD = C$ has, up to a permutation of rows, the
matrix given in Table~\ref{EPims} as its only solution. (We have already
chosen the ordering with some foresight such that a later reordering of
rows or columns is unnecessary.)

\begin{table}
\[ 
\begin{array}{c|ccccc} \hline
 & \Phi_{3a} & \Phi_{1a} & \Phi_{1b} & \Phi_{1c} & \Phi_{2a} \\ 
 & 3a & 1a & 1b & 1c & 2a \\ \hline
3_1 & 1 & . & . & . & . \rule[ .pt]{.pt}{ 13pt} \\
1_1 & . & 1 & . & . & . \\
5_1 & 1 & 1 & 1 & . & . \\
4_1 & 1 & . & . & 1 & . \\
1_2 & . & . & . & 1 & . \\
5_2 & 1 & . & 1 & 1 & . \\
3_2 & . & . & 1 & . & 1 \\
4_2 & . & 1 & 1 & . & 1 
\rule[- 7pt]{0pt}{ 5pt} \\ \hline
\end{array}
\hspace*{5em}
C := \left( \begin{array}{ccccc}
  4&  1&  2&  2&  . \\
  1&  3&  2&  .&  1 \\
  2&  2&  4&  1&  2 \\
  2&  .&  1&  3&  . \\
  .&  1&  2&  .&  2  \\
\end{array} \right) 
\]
\caption{\label{EPims} Decomposition matrix and Cartan matrix of~$E$}
\end{table}

Hence it remains to show how the rows of the matrix in 
Table~\ref{EPims} correspond to rows of Table~\ref{PS1}. 
The regular character of $E$ shows that 
\[ \begin{array}{ll}
\{1_1,1_2\}\leftrightarrow\{\chi_{10},\chi_{37}\}, &
\{3_1,3_2\}\leftrightarrow\{\chi_{8},\chi_{49}\}, \\
\{4_1,4_2\}\leftrightarrow\{\chi_{33},\chi_{50}\}, &
\{5_1,5_2\}\leftrightarrow\{\chi_{32},\chi_{43}\}, \\
\end{array} \]
while $\chi_{19}$ is not met. 
The basic set of Table~\ref{PS1} shows in particular that $\chi'_8$ is 
an irreducible Brauer character, and the associated projective indecomposable
character $\Phi_1$ is contained three times as a summand in the 
projective character $\Psi$. Thus~$\Phi_{3a}$ is the Fitting 
correspondent of~$\Phi_1$, and part of the matching is 
\[ 3_1\leftrightarrow\chi_{8} \quad\text{and}\quad 
   3_2\leftrightarrow\chi_{49} .\]
This yields $8$ remaining cases.
To be a valid matching, the projective indecomposable characters
of $E$ have to decompose into the basic set of projective characters
in Table~\ref{PS1} with integral coefficients. It is
easily checked using \textsf{GAP}
that this necessary condition is only fulfilled for
\[ \begin{array}{llllll}
1_1\leftrightarrow\chi_{10}, &
1_2\leftrightarrow\chi_{37}, &
4_1\leftrightarrow\chi_{33}, &
4_2\leftrightarrow\chi_{50}, &
5_1\leftrightarrow\chi_{32}, &
5_2\leftrightarrow\chi_{43}. \\
\end{array} \]
This gives five of the seven projective indecomposable characters 
for block~$B_1$, and together with $\Psi_3$ and $\Psi_6$
we obtain the basic set of projective characters displayed in Table~\ref{PS2}.

\begin{table}
\[
\begin{array}{r|r|rrrrrrr|r} \hline
\chi & \chi(1) & \Phi_1 & \Phi_2 & \Psi_3 & \Phi_4 & \Phi_5 & \Psi_6 & \Phi_7 & \Psi' \\ \hline
{\bf 8} &   8\,910 &  1 &   . &   . &   . &   . &   . &   . & . 
 \rule[ .pt]{.pt}{ 13pt} \\
{\bf 10} &  16\,929 &  . &   1 &   . &   . &   . &   . &   . & . \\
{\bf 19} & 270\,864 &  . &   . &   1 &   . &   . &   . &   . & 1 \\
{\bf 32} &1\,185\,030 &  1 &   1 &   . &   1 &   . &   1 &   . & 1 \\
{\bf 33} &1\,354\,320 &  1 &   . &   . &   . &   1 &   2 &   . & . \\
{\bf 37} &1\,575\,936 &  . &   . &   1 &   . &   1 &   1 &   . & 1 \\
43&2\,784\,375 &  1 &   . &   1 &   1 &   1 &   3 &   . & 2 \\
{\bf 49} &4\,561\,920 &  . &   . &   2 &   1 &   . &   6 &   1 & 2 \\
50&4\,809\,375 &  . &   1 &   3 &   1 &   . &   5 &   1 & 3
\rule[- 7pt]{0pt}{ 5pt} \\ \hline
\end{array}
\]
\caption{\label{PS2} A second basic set of projective characters for~$B_1$}
\end{table}

\medskip
The projective character $\Psi'$ given in the last column of
Table~\ref{PS2} is obtained as $(\chi_2 \cdot \chi_{30})\cdot 1_{B_1}$.
We have \[\Psi' = \Psi_3 + \Phi_4 - \Phi_7.\] Since $\Phi_7$ and
$\Phi_4$ are projective indecomposable, this relation implies that
\[\Psi_3' := \Psi_3 - \Phi_7\] is a projective character. We obtain the
new basic set of projective characters displayed in Table~\ref{PS3}.
>From this table it is evident that in order to obtain the projective
indecomposable character $\Phi_3$ contained in $\Psi_3'$, we may
subtract $\Phi_7$ at most once from $\Psi_3'$. Likewise, we consider the
decomposition
\[ \Phi_6 = \Psi_6 - a \cdot \Phi_4 -b\cdot \Phi_5 - c\cdot\Phi_7, 
\quad\text{where }0\le a, b \le 1\text{ and }0 \le c \le 5 - a .\]

This leaves 44 possibilities for the decomposition matrix of $B_1$.
We postpone the proof which of these possibilities holds until 
after our treatment of the principal block, as with little 
effort we can eliminate three quarters of these possibilities 
in parallel to determining the Brauer characters of $B_0$.

\begin{table}
\[
\begin{array}{r|r|rrrrrrr} \hline
\chi & \chi(1) & \Phi_1 & \Phi_2 & \Psi_3' & \Phi_4 & \Phi_5 & \Psi_6 & \Phi_7 \\ \hline
{\bf 8} &   8\,910 &  1 &   . &   . &   . &   . &   . &   .  
  \rule[ .pt]{.pt}{ 13pt} \\
{\bf 10} &  16\,929 &  . &   1 &   . &   . &   . &   . &   .  \\
{\bf 19} & 270\,864 &  . &   . &   1 &   . &   . &   . &   .  \\
{\bf 32} &1\,185\,030 &  1 &   1 &   . &   1 &   . &   1 &   .  \\
{\bf 33} &1\,354\,320 &  1 &   . &   . &   . &   1 &   2 &   .  \\
{\bf 37} &1\,575\,936 &  . &   . &   1 &   . &   1 &   1 &   .  \\
43&2\,784\,375 &  1 &   . &   1 &   1 &   1 &   3 &   .  \\
{\bf 49} &4\,561\,920 &  . &   . &   1 &   1 &   . &   6 &   1  \\
50&4\,809\,375 &  . &   1 &   2 &   1 &   . &   5 &   1 
\rule[- 7pt]{0pt}{ 5pt} \\ \hline
\end{array}
\]
\caption{\label{PS3} A third basic set of projective characters for~$B_1$}
\end{table}

\subsection{Proof for the principal block $B_0$}
\label{subsec:B0mod3}
As in our treatment of $G$ in characteristic 2, we begin by constructing
several small dimensional representations. In \cite{web:ATLAS} we
readily find $133_1$ and $133_2$, which are conjugate under the outer
automorphism of $G$, and $760$. Furthermore, restricting the smallest
non-trivial irreducible representation of the sporadic Baby Monster
group $\mathsf{B}$ also available from \cite{web:ATLAS} gives the
composition factors
\[
4\,371\smash\downarrow_G = 1 + 133_1 + 133_2 + 760 + 3\,344,\]
all of which are liftable.
As the tensor products $133_1 \otimes 133_1$ and $133_1 \otimes 133_2$ allow
a direct treatment with the MeatAxe, we also chop these modules obtaining
\begin{eqnarray*}
 133_1 \otimes 133_1 &=& 1 + 8\,778_1 + 8\,910\\
 133_1 \otimes 133_2 &=& 760 + 16\,929
\end{eqnarray*}
in the process. Choosing basic sets of Brauer characters for
$B_0$ and $B_1$ as indicated in bold face in Tables~\ref{PS3} and
\ref{tab:decB0mod3}, respectively, we see that $8\,910$ and $16\,929$
afford two of the irreducible Brauer characters of $B_1$ already known.
Hence, we readily obtain the Brauer character of $8\,778_1$, and
applying the outer automorphism to $8\,778_1$ gives $8\,778_2$.

\medskip
To determine the remaining Brauer characters, we employ condensation
again. As condensation subgroups we choose $K_1\cong 2^3.2^2.2^6$ of
order $2\,048$, which is the largest normal 2-subgroup of the ninth
maximal subgroup $M_9\cong 2^3.2^2.2^6.(3\times L_3(2))$ of $G$, and
$K_2\cong 5^2.5.5^2$ of order $3\,125$, which is the largest normal
5-subgroup of the tenth maximal subgroup $M_{10}\cong 5^2.5.5^2.4A_5$
of $G$ (see \cite{ATLAS}). A straight line programme to find $M_{10}$
is available in \cite{web:ATLAS}. Following the construction given in
\cite[Section 3, p.\ 365]{NortonWilson} we can similarly procure a
straight line programme for the ninth maximal subgroup by computing its
intersection with the maximal subgroup $A_{12}$. The normal subgroups
in question are then easily found using the facilities to deal with
permutation groups provided by \textsf{GAP}.

The reason for choosing two condensation subgroups becomes evident, if
we calculate the dimensions of the condensed simple modules known via
their Brauer characters: we see that the corresponding trace idempotent
$e_1$ annihilates the pairs $(133_1, 133_2)$ and $(8\,778_1, 8\,778_2)$,
whereas the trace idempotent $e_2$ of the second condensation subgroup
$K_2$ annihilates $760$ and $3\,344$. Condensing every module with
both idempotents allows us to counteract the blind spots produced by
condensing with just either of the two.

\medskip
By \cite[Theorem 2.7]{tackling} condensing the non-identity
representatives of the double cosets in $M_9 \backslash G /M_9$ and
$M_{10}\backslash G /M_{10}$ together with generators for the maximal
subgroups gives generating sets for $e_1FGe_1$ and $e_2FGe_2$,
consisting of 387 and 643 elements, respectively. Words for the double
coset representatives can be computed using \textsf{orb}, realising the
permutation representations of $G$ on the 264\,515?\,625 cosets in $M_9
\backslash G$, and the 364\,041\,216 cosets in $M_{10}\backslash G$, as
$G$-orbits of vectors in suitable matrix representations of $G$; again
we spare the details of the actual computations.

As working with 1030 generators simultaneously is cumbersome to
impossible, we adopt the strategy to consider the algebras $\Ccal_1 \le
e_1FGe_1$ and $\Ccal_2 \le e_2FGe_2$, each generated by the condensed
elements corresponding to the following six elements
\begin{equation}\label{eq:sixgens} 
 a, \: b, \: ab, \: ba, \: a^2, \: aba ,
\end{equation}
where $a$ and $b$ are standard generators of $G$. We apply
the condensation of tensor products technique described in
\cite{art:newtensorcondense, LuxWie} to the modules listed in
Table~\ref{tab:condresmod3}, where the simple module $12\,264_1$ which
we also use will be constructed in the process, and chop the condensed
modules using the six group elements of \eqref{eq:sixgens}.

\medskip
For any $FG$-module $M$ the idea now is to consider a composition
series of $Me_i\smash\downarrow_{\Ccal_i}$ for $i=1,2$, which can be
obtained with the MeatAxe. Let $\Bcal_i$ be a basis of $Me_i$ adapted
to the composition series computed, i.e.\ with respect to $\Bcal_i$ the
action of elements of $\Ccal_i$ on $Me_i\smash\downarrow_{\Ccal_i}$ is
given by lower block-triangular matrices, for which the block-diagonal
gives the representations on the composition factors. We can prove that
every composition factor of $Me_i$ restricts irreducibly to $\Ccal_i$
by checking if the lower block-diagonal structure is preserved by
all generators of $e_iFGe_i$ with respect to the basis $\Bcal_i$.
Moreover, we check that two composition factors are isomorphic as
$e_iFGe_i$-modules if they are as $\Ccal_i$-modules. This is applied to
the tensor products $760 \otimes 8\,778_1$ and $3\,344 \otimes 3\,344$.
Thus all $\Ccal_1$- and $\Ccal_2$-modules appearing as composition
factors in the modules can be uniquely extended to $e_1FGe_1$- and
$e_2FGe_2$-modules, respectively.

\medskip
For the identification which simple $e_1FGe_1$-module corresponds
to which simple $e_2FGe_2$-module, we now adopt the improved matching 
method detailed in \cite{matching2}. 

\medskip
Using our basic sets of $B_0$ and $B_1$, together with the Brauer
characters of $B_2$ and the characters of defect zero, it is plain
to see which composition factors outside of $B_0$ and $B_1$ appear
in the tensor products $760 \otimes 8\,778_1$ and $3\,344 \otimes
3\,344$. The dimensions of the their condensed modules are
also readily determined. These dimensions are almost sufficient to
identify these modules, and thus to omit them immediately from the
condensation results, except that there is a constituent of dimension
$96$ for the condensation subgroup $K_2$ being the condensed module of
the irreducible Brauer character $406\,296$ in the block $B_2$ of defect
$1$. Running the matching algorithm shows that this constituent matches
with a constituent of dimension $243$ for the condensation subgroup
$K_1$, hence can be uniquely identified. 
Moreover, the decompositions of both tensor products
into our basic sets of $B_0$ and $B_1$ (see Table~\ref{tab:decBSmod3})
show that the character $\chi'_{49}$ of the basic set of Table~\ref{PS3} does
not occur in either decomposition. Hence by the results obtained in
\ref{subsec:B1mod3part1}, the irreducible Brauer character corresponding
to the projective indecomposable $\Phi_7$ is not a constituent of either
of these tensor products. Thus their restrictions to $B_0 \oplus B_1$
contain at most $26$ distinct constituents.

The result of the matching algorithm run on the components in the
condensed blocks of $B_0$ and $B_1$ is given in the second and third
columns of Table~\ref{tab:condresmod3}. More precisely, we indeed
find the matching between the constituents as shown, where only
$k1_1$ and $k_9$ for the condensation subgroup $K_1$, and $c3_{1/2}$,
$c8_{1/2}$ and $c16_{1/2}$ for the condensation subgroup $K_2$ remain
`unmatched'. This means that either we have missed existing matchings,
or there is none because of non-faithful condensation. Since we have
already constructed the simple modules $133_{1/2}$, $760$, $3\,344$,
$8\,778_{1/2}$, and will do so for the simple modules $12\,264_{1/2}$
later on, we can identify all of the `unmatched' constituents as
condensed simple modules, and thus show by using their Brauer characters
that there indeed is no matching for them.

Thus we count 26 non-isomorphic simple modules occurring as
`matched composition factors' in the condensation results. Hence, all
simple modules of the principal block $B_0$, and all but one simple module
of $B_1$ arise as composition factors in the two tensor products.
Furthermore, we have a bijection between these simple $FG$-modules and
the above `matched composition factors'.

\begin{table}
 \[
 \begin{array}{r|r|r|rrrrrrrrrrrr}
  \hline
& K_1 & K_2 & 
  \begin{turn}{90}
   $133_1 \otimes 133_1$
  \end{turn} &
  \begin{turn}{90}
   $133_1 \otimes 133_2$
  \end{turn} &
  \begin{turn}{90}
   $133_1 \otimes 3\,344$
  \end{turn} &
  \begin{turn}{90}
   $133_1 \otimes 760$
  \end{turn} &
  \begin{turn}{90}
   $760 \otimes 760$
  \end{turn} &
  \begin{turn}{90}
   $1_{A_{12}}^G$
  \end{turn} &
  \begin{turn}{90}
   $760 \otimes 3\,344$
  \end{turn} &
  \begin{turn}{90}
   $1_{2.\mathsf{HS}.2}^G$
  \end{turn} &
  \begin{turn}{90}
   $133_1 \otimes 8\,778_1$
  \end{turn} &
  \begin{turn}{90}
   $3\,344 \otimes 3\,344$
  \end{turn} &
  \begin{turn}{90}
   $133_1 \otimes 12\,264_1$
  \end{turn} &
  \begin{turn}{90}
   $760 \otimes 8\,778_1$
  \end{turn}\\ \hline
760 & k1_1 & - & . & 1 & . & 1 & . & 2 & 2 & 3 & 5 & 16 & 1 & 3\\
1 & k1_2 & c1 & 1 & . & . & 1 & 1 & 2 & 1 & 4 & 1 & 6 & . & 3\\
133_2 & - & c3_1 & . & . & . & 1 & 1 & 2 & . & 2 & . & 4 & . & 1\\
133_1 & - & c3_2 & . & . & . & . & 1 & 2 & . & 2 & 2 & 4 & . & 1\\
9\,139 & k8 & c7 & . & . & . & . & 2 & 3 & 1 & 7 & 2 & 13 & 1 & 3\\
3\,344 & k9 & - & . & . & 2 & 1 & 2 & 5 & 4 & 6 & 2 & 15 & . & 3\\
12\,264_1 & - & c8_1 & . & . & . & 1 & . & 1 & . & 2 & 3 & 10 & 1 & 2\\
12\,264_2 & - & c8_2 & . & . & . & 1 & . & 1 & . & 2 & 3 & 10 & 1 & 2\\
31\,768_1 & k26_1 & c8_3 & . & . & 1 & . & . & 2 & 2 & 3 & 2 & 9 & 1 & 3\\
31\,768_2 & k26_2 & c8_4 & . & . & . & 2 & . & 2 & 2 & 3 & 1 & 9 & . & 4\\
8\,778_1 & - & c16_1 & 1 & . & . & . & 1 & . & . & . & 4 & 4 & . & .\\
8\,778_2 & - & c16_2 & . & . & . & 1 & 1 & . & . & . & . & 4 & . & .\\
137\,236 & k50 & c40 & . & . & . & . & . & . & 1 & 2 & . & 4 & 1 & 2\\
147\,061 & k89 & c33 & . & . & . & . & . & . & . & 1 & . & 2 & 1 & 2\\
339\,702_2 & k138_1 & c96_2 & . & . & . & . & . & . & . & . & 2 & 2 & . & .\\
339\,702_1 & k138_2 & c96_1 & . & . & . & . & . & . & . & . & . & 2 & . & .\\
255\,037 & k173 & c89 & . & . & . & . & 1 & 2 & . & 3 & 1 & 6 & . & .\\
496\,924_1 & k260_1 & c168_1 & . & . & . & . & . & . & . & . & . & . & 1 & 1\\
496\,924_2 & k260_2 & c168_2 & . & . & . & . & . & . & . & . & . & . & . & 1\\
783\,696 & k322 & c272 & . & . & . & . & . & . & . & . & . & 1 & 1 & 1\\
  \hline
40\,338 & k18 & c6 & . & . & . & . & . & . & . & . & . & . & . & 1\\
8\,910 & k27_1 & c18 & 1 & . & . & . & 1 & 1 & . & 1 & . & 2 & . & 1\\
16\,929 & k27_2 & c9 & . & 1 & . & . & . & 1 & . & 1 & . & 2 & . & 1\\
270\,864 & k105 & c96 & . & . & . & . & 1 & . & 2 & . & . & 2 & . & 1\\
1\,159\,191 & k618 & c387 & . & . & . & . & . & . & . & . & . & 1 & . & .\\
1\,305\,072 & k702 & c384 & . & . & . & . & . & . & 1 & . & . & 1 & . & 2\\
  \hline
 \end{array}
 \]
 \caption{Condensation results for $p=3$ for $\mathsf{HN}$}
 \label{tab:condresmod3}
\end{table}

\begin{table}
 \[
 \begin{array}{r|r|rrrrrrrrrrrr|r}
\hline \chi & \chi(1) &
  \begin{turn}{90}
   $133_1 \otimes 133_1$
  \end{turn} &
  \begin{turn}{90}
   $133_1 \otimes 133_2$
  \end{turn} &
  \begin{turn}{90}
   $133_1 \otimes 3\,344$
  \end{turn} &
  \begin{turn}{90}
   $133_1 \otimes 760$
  \end{turn} &
  \begin{turn}{90}
   $760 \otimes 760$
  \end{turn} &
  \begin{turn}{90}
   $1_{A_{12}}^G$
  \end{turn} &
  \begin{turn}{90}
   $760 \otimes 3\,344$
  \end{turn} &
  \begin{turn}{90}
   $1_{2.\mathsf{HS}.2}^G$
  \end{turn} &
  \begin{turn}{90}
   $133_1 \otimes 8\,778_1$
  \end{turn} &
  \begin{turn}{90}
   $3\,344 \otimes 3\,344$
  \end{turn} &
  \begin{turn}{90}
   $133_1 \otimes 12\,264_1$
  \end{turn} &
  \begin{turn}{90}
   $760 \otimes 8\,778_1$
  \end{turn} &
  \begin{turn}{90}
   $3\,344 \otimes 9\,139_1$
  \end{turn}\\ 
\hline
  1 &       1 &1& .& .&  .& 1&  1& .&  1&  .&  1&  1&  1 &  2\\
  2 &     133 &.& .& .&  .& .&  1& .&  .&  1&  .& -2& -1 & -1\\
  3 &     133 &.& .& .&  1& .&  1& .&  .& -1&  .& -2& -1 & -1\\
  4 &     760 &.& 1& .&  .& .&  1& 1&  .&  .&  .& -1&  . & -2\\
  5 &    3344 &.& .& 1&  .& 1&  .& 1&  .& -1& -2& -1& -2 & -2\\
  6 &    8778 &1& .& .& -1& 1& -1& .& -1&  1& -1&  1&  . &  1\\
  7 &    8778 &.& .& .&  1& 1&  .& .&  .&  .&  1& -2& -1 & -3\\
  9 &    9405 &.& .& .&  .& 1&  1& .&  2&  1&  4&  2&  2 &  3\\
 11 &   35112 &.& .& 1& -1& .&  1& 1&  .&  1&  3&  3&  2 &  3\\
 12 &   35112 &.& .& .&  2& .&  2& 1&  1&  1&  5& -1&  2 & -1\\
 13 &   65835 &.& .& .&  1& .&  1& .&  1&  1&  2& -3& -1 & -4\\
 17 &  214016 &.& .& .&  .& .&  .& 1&  1&  .&  1&  .&  1 &  .\\
 18 &  267520 &.& .& .&  .& 1&  2& .&  2&  1&  4&  .&  . &  .\\
 21 &  374528 &.& .& .&  .& .&  .& .&  .&  .&  2&  1&  . &  2\\
 22 &  374528 &.& .& .&  .& .&  .& .&  .&  2&  2&  1&  . &  2\\
 24 &  653125 &.& .& .&  .& .&  .& .&  1&  .&  2&  .&  . &  3\\
 25 &  656250 &.& .& .&  .& .&  .& .&  .&  .&  .&  1&  1 &  1\\
 26 &  656250 &.& .& .&  .& .&  .& .&  .&  .&  .&  1&  1 &  1\\
 29 & 1053360 &.& .& .&  .& .&  .& .&  .&  .&  1&  1&  1 &  1\\
 35 & 1361920 &.& .& .&  .& .&  .& .&  .&  .&  .& -1&  . &  .\\
   \hline                                                   
  8 &    8910 &1& .& .&  .& 1&  1& .&  1&  .&  1&  .&  . &  .\\
 10 &   16929 &.& 1& .&  .& .&  1& .&  1&  .&  1&  .&  1 &  .\\
 19 &  270864 &.& .& .&  .& 1&  .& 1&  .&  .&  1&  .&  . &  1\\
 32 & 1185030 &.& .& .&  .& .&  .& .&  .&  .&  1&  .&  . &  1\\
 33 & 1354320 &.& .& .&  .& .&  .& .&  .&  .&  .&  .&  1 &  .\\
 37 & 1575936 &.& .& .&  .& .&  .& 1&  .&  .&  1&  .&  1 &  1\\
 49 & 4561920 &.& .& .&  .& .&  .& .&  .&  .&  .&  .&  . &  1\\
   \hline
\end{array}
\]
\caption{Decompositions into basic sets for $B_0$ and $B_1$}
\label{tab:decBSmod3}
\end{table}

\medskip
The strategy now is as follows: We proceed through the columns of
Tables~\ref{tab:condresmod3} and \ref{tab:decBSmod3}, 
and compare the multiplicities of the
`matched composition factors' with the decomposition of the Brauer
characters of the tensor products in question into the basic sets of
$B_0$ and $B_1$. This way we collect successively the information how
the basic set characters decompose into irreducible Brauer characters.
Using this information, it is straightforward to determine the Brauer
characters which correspond to the `matched composition factors'. We
therefore omit the elementary calculations and only state the individual
results:

\medskip
$\blacktriangleright$ 
Condensing the simple modules, which are explicitly known
as matrix representations, first reveals 
\[ 1e_1\mapsto k1_2,\: 760e_1\mapsto k1_1,\: 3\,344 e_1\mapsto k9,\:
   133_1 e_2\mapsto c3_2,\: 133_2e_2\mapsto c3_1 .\]

\medskip
$\blacktriangleright 133_1\otimes 133_1$ and $133_1\otimes 133_2$: 
These are condensed to match the composition factors,
which we have already determined but not computed as explicit
matrix representations:
\[ 8\,778_1e_2\mapsto c16_1,\: 8\,778_2e_2\mapsto c16_2,\:
   8\,910e_1\mapsto k27_1,\: 16\,929e_1\mapsto k27_2 .\]

\medskip
$\blacktriangleright 133_1 \otimes 3\,344$: 
This yields the Brauer characters $31\,768_1$ and its conjugate
$31\,768_2$:
\[ 31\,768_1e_1\mapsto k26_1,\: 31\,768_2e_1\mapsto k26_2 .\]

\medskip
$\blacktriangleright 133_1 \otimes 760$: 
The only basic set character needed to express the Brauer
character of $133_1 \otimes 760$ whose modular constituents 
are not yet known is $\chi_{13}$, of degree $65\,835$. Hence we obtain
\begin{equation}\label{eq:sum12264}
 \chi_{13}' = 1 + 760 + 8\,778_1 + c8_1 + c8_2 + 31\,768_1,
\end{equation}
where we denote the unknown Brauer characters in this decomposition
by their corresponding condensed modules. To determine these Brauer
characters, we analyse the condensed tensor product with the MeatAxe.
Its socle series is revealed to be
\begin{equation}
 \boxed{
 \begin{matrix}
  c8_4 \\
  c8_2 \\
  c16_2 \\
  c8_1 \oplus c1\\
  c8_4 \oplus c3_1
 \end{matrix}}\:.
 \label{eq:soc133x760}
\end{equation}
Analogously to our approach in
Section~\ref{subsec:B0mod2}, we construct a basis for the uncondensed
submodule of Loewy length two, with head $c8_1$ and socle $c8_4 \oplus
c3_1$. With it we construct a representation of the module condensing
to $c8_1$, giving $12\,264_1$. As $133_1 \otimes 760$ is self-dual, but
the contragredient of $12\,264_1$ gives a new simple module $12\,264_2$,
which can be checked with the representation,
we can derive from the socle series in \eqref{eq:soc133x760} that
$$ 12\,264_1e_2\mapsto c8_1,\: 12\,264_2e_2\mapsto c8_2 .$$

>From \eqref{eq:sum12264} we explicitly obtain the Brauer character
of the sum $12\,264_1 + 12\,264_2$. 
Let $\sigma$ denote the trivial extension of the nontrivial Galois
automorphism of $\mathbb{Q}(\sqrt{5})$ to the algebraic number field
containing all values of the ordinary characters.
The character table of $G$ (see \cite{ATLAS})
reveals that the outer automorphism acts on the ordinary characters 
as the product of $\sigma$ and complex conjugation.
Now, since $12\,264_2$ is the
contragredient of $12\,264_1$, and both are interchanged by the action
of the outer automorphism (which can again be checked with the
representations obtained and a straightline programme for the action of the
outer automorphism available in \cite{web:ATLAS}), 
their characters are invariant under the
action of the Galois group of $\mathbb{Q}(\sqrt{5})$. 
It therefore follows that
$12\,264_1$ and $12\,264_2$ are rational on all classes except the
pairs $19A/B$ and $40A/B$ (note that irrationalities on $35A/B$ only 
appear as values of defect zero characters). 
Since we have constructed representations
for $12\,264_1$, and straight line programs for conjugacy class
representatives of $G$ are available in \cite{web:ATLAS}, we can compute
representatives for $19A$ and $40A$, and the Brauer character value
of the respective representative. We find that $12\,264_1$ assumes
the value $b19$ on $19A$ and $i10$ on $40A$, using the notation in
\cite{ATLAS}. The values on all classes with rational values are simply
obtained by dividing all the values of the character sum $12\,264_1 +
12\,264_2$ by 2. Thus, this construction yields the Brauer characters
$12\,264_1$ and $12\,264_2$.

\medskip
$\blacktriangleright 760\otimes 760$ and $1_{A_{12}}^G$: Considering both 
columns of Tables~\ref{tab:condresmod3} and \ref{tab:decBSmod3} gives
the equations
\begin{eqnarray*}
 \chi_9' + \chi_{18}' &=& 133_1 + 133_2 + 3\,344 + 2\times k8 + k173,\\
 \chi_9' + 2 \chi_{18}' &=& 133_1 + 133_2 + 2\times 3\,344 + 3\times
 k8 + 2\times k173.
\end{eqnarray*}
This allows us to derive the constituents of $\chi_9'$ and 
$\chi_{18}'$ as
\begin{eqnarray*}
 \chi_9' &=& 9\,139 + 133_1 + 133_2,\\
 \chi_{18}' &=& 255\,037 + 9\,139 + 3\,344,
\end{eqnarray*}
and considering the component of $760\otimes 760$ in $B_1$
we recover the Brauer character $270\,864$, hence we see that 
$$ 9\,139\,e_1\mapsto k8,\: 255\,037\, e_1\mapsto k173,\:
   270\,864\,e_1\mapsto k105 .$$

\medskip
$\blacktriangleright 760\otimes 3\,344$ and 
$1_{2.\mathsf{HS}.2}^G:$ From $1_{2.\mathsf{HS}.2}^G$ we obtain the relation
\[
\begin{split}
 k50+k89 = \chi_{24}' &-1 -2\times 9\,139 - 2\times 3\,344 
-12\,264_1 -12\,264_2\\
 &-31\,768_1 - 31\,768_2 -255\,037.
\end{split}
\]
In particular we derive from this decomposition that both $k89$ and
$k50$ belong to the principal block.

Now considering $760\otimes 3\,344$, from the decomposition into the
basic sets in Table~\ref{tab:decBSmod3}
we infer that, since $k50$ lies in $B_0$, the other occurring
unknown simple module $k702$ lies in $B_1$. With the help of the basic
sets we may separate the character of the tensor product restricted to
$B_0$ and $B_1$ into its block components, and obtain via
\begin{eqnarray*}
 k50 &=& \chi_{17}' -1 -760- 3\,344- 31\,768_1 -31\,768_2 -9\,139,\\
 k702 &=& \chi_{37}' - 270\,864
\end{eqnarray*}
the Brauer characters $137\,236$ of $B_0$ and $1\,305\,072$ of $B_1$.

Reconsidering the above relation derived from $1_{2.\mathsf{HS}.2}^G$
yields the Brauer character $147\,061$ of $B_0$, thus
$$ 137\,236\,e_1\mapsto k50,\: 147\,061\mapsto k89,\:
   1\,305\,072\,e_1\mapsto k702 .$$

\medskip
$\blacktriangleright 133_1 \otimes 8\,778_1$:
Here we obtain the decomposition
\[ \chi_{22}' = 2\times 760 + 12\,264_1 + 12\,264_2 + 8\,778_1 + k138_1.
\]
Hence we get an irreducible Brauer character of degree $339\,702$, which we
call $339\,702_2$ as it appears in the second basic set character of degree
$374\,528$. Applying the outer automorphisms gives $339\,702_1$, hence
$$ 339\,702_2\,e_1\mapsto k138_1,\: 339\,702_1\,e_1\mapsto k138_2 .$$

\medskip
$\blacktriangleright 3\,344 \otimes 3\,344$:
Again separating the character into its block components with our basic
sets, we see that this tensor product yields one yet unknown simple
module of $B_0$ and likewise one yet unknown simple module of $B_1$.
Thus we obtain a Brauer character $783\,696$ of $B_0$, and another
Brauer character $1\,159\,191$ of $B_1$. Computing the dimensions of the
corresponding simple condensed modules with these characters, we see
that $$ 783\,696\,e_1\mapsto k322,\: 1\,159\,191\,e_1\mapsto k618 .$$

\medskip
$\blacktriangleright 133_1 \otimes 12\,264_1$:
This tensor product gives the relation
\[ \chi_{25}' + \chi_{26}' - \chi_{35}'= k260_1 + 1 + 12\,264_1 + 147\,061.\]
Therefore we have a new Brauer character $496\,924_1$. 
>From the latter we obtain its contragredient $496\,924_2$, hence
$$ 496\,924_1\,e_1\mapsto k260_1,\: 496\,924_2\,e_1\mapsto k260_2 .$$
Thus we have determined all irreducible Brauer
characters of $B_0$ at this stage, and we give the 3-modular
decomposition matrix of this Block in Table~\ref{tab:decB0mod3}.

\medskip
$\blacktriangleright 760 \otimes 8\,778_1$:
Finally, considering the last column of 
Table~\ref{tab:condresmod3} we obtain another
Brauer character of $B_1$: we derive the decomposition
\[ \chi_{33}' = k18 + 1\,305\,072 + 8\,910,\]
which immediately yields the Brauer character $40\,338$, hence
$$ 40\,338\,e_1\mapsto k18 .$$

\begin{landscape}
\begin{table}
 {\tiny
 \[
 \begin{array}{r|r|rrrrrrrrrrrrrrrrrrrr}
  \hline
\chi & \chi(1) & \Phi_1 & \Phi_2 & \Phi_3 & \Phi_4 & \Phi_5 & \Phi_6 &
\Phi_7 & \Phi_8 & \Phi_9 & \Phi_{10} & \Phi_{11} & \Phi_{12} &
\Phi_{13} & \Phi_{14} & \Phi_{15} & \Phi_{16} & \Phi_{17} &
\Phi_{18} & \Phi_{19} & \Phi_{20} \\
\hline
{\bf 1}& 1 & 1& .& .& .& .& .& .& .& .& .& .& .& .& .& .& .& .& .& .& . \\ 
{\bf 2}& 133 & .& 1& .& .& .& .& .& .& .& .& .& .& .& .& .& .& .& .& .& . \\ 
{\bf 3}& 133 & .& .& 1& .& .& .& .& .& .& .& .& .& .& .& .& .& .& .& .& . \\ 
{\bf 4}& 760 & .& .& .& 1& .& .& .& .& .& .& .& .& .& .& .& .& .& .& .& . \\ 
{\bf 5}& 3\,344 & .& .& .& .& 1& .& .& .& .& .& .& .& .& .& .& .& .& .& .& . \\ 
{\bf 6}&8\,778 & .& .& .& .& .& 1& .& .& .& .& .& .& .& .& .& .& .& .& .& . \\ 
{\bf 7}&8\,778 & .& .& .& .& .& .& 1& .& .& .& .& .& .& .& .& .& .& .& .& . \\ 
{\bf 9}&9\,405 & .& 1& 1& .& .& .& .& 1& .& .& .& .& .& .& .& .& .& .& .& . \\ 
{\bf 11}&35\,112&.&.&.& .& 1& .& .& .& .& .& 1& .& .& .& .& .& .& .& .& . \\ 
{\bf 12}&35\,112&.&.&.& .& 1& .& .& .& .& .& .& 1& .& .& .& .& .& .& .& . \\ 
{\bf 13}&65\,835&1&.&.&1& .& 1& .& .& 1& 1& 1& .& .& .& .& .& .& .& .& . \\ 
14& 65\,835 & 1& .& .& 1& .& .& 1& .& 1& 1& .& 1& .& .& .& .& .& .& .& . \\ 
{\bf 17}&214\,016 &1&.&.&1& 1& .& .& 1& .& .& 1& 1& 1& .& .& .& .& .& .& . \\ 
{\bf 18}&267\,520 &.&.&.& .& 1& .& .& 1& .& .& .& .& .& .& 1& .& .& .& .& . \\ 
20& 365\,750 & 1& 1& 1& 1& 1& .& .& 2& 1& 1& 1& 1& .& .& 1& .& .& .& .& . \\ 
{\bf 21}&374\,528 &.&.&.&2&.& .& 1& .& 1& 1& .& .& .& .& .& 1& .& .& .& . \\ 
{\bf 22}&374\,528 &.&.&.&2&.& 1& .& .& 1& 1& .& .& .& .& .& .& 1& .& .& . \\ 
{\bf 24}&653\,125 &1&.&.&1&2& .& .& 2& 1& 1& 1& 1& 1& 1& 1& .& .& .& .& . \\ 
{\bf 25}&656\,250 &1&.&.& .& .& .& .& .& 1& .& .& .& .& 1& .& .& .& 1& .& . \\ 
{\bf 26}&656\,250 &1&.&.&.&.&.& .& .& .& 1& .& .& .& 1& .& .& .& .& 1& . \\ 
27& 718\,200 & .& .& 1& 2& 2& .& 1& 2& 1& 1& 1& 1& .& .& 1& 1& .& .& .& . \\ 
28& 718\,200 & .& 1& .& 2& 2& 1& .& 2& 1& 1& 1& 1& .& .& 1& .& 1& .& .& . \\ 
{\bf 29}&1\,053\,360 &.&.&.&3&.&1&1& .& 2& 2& 1& 1& 1& .& .& .& .& .& .& 1 \\ 
34& 1\,361\,920 & .& .&.&3& 3& .& .& 2& 2& 2& 1& 1& 1& 1& 1& 1& 1& .& .& . \\ 
{\bf 35}&1\,361\,920 &.&.&.&2&1&.& .& 1& 1& 1& .& .& .& 1& .& 1& 1& .& 1& . \\ 
36& 1\,361\,920 & .& .& .& 2& 1& .& .& 1& 1& 1& .& .& .& 1& .& 1& 1& 1& .& . \\ 
40& 2\,375\,000 & 2& 1& 1& 3& 3& .& .& 3& 2& 2& 1& 1& .& 2& 1& 1& 1& 1& 1& . \\ 
41& 2\,407\,680 & 2& 1& 1& 1& 1& .& .& 2& 2& 2& 2& 2& 1& 2& .& .& .& 1& 1& 1 \\ 
42& 2\,661\,120 & 2& 1& 1& 4& 2& .& .& 4& 3& 3& 1& 1& .& 2& 2& 1& 1& 1& 1& . \\ 
45& 3\,200\,000 & 1& 1& 1& 3& 2& .& .& 2& 3& 3& 1& 1& 1& 3& .& 1& 1& 1& 1& 1 \\ 
46& 3\,424\,256 & 2& 1& 1& 5& 2& 1& 1& 3& 4& 4& 2& 2& 1& 2& 1& 1& 1& 1& 1& 1 \\ 
48& 4\,156\,250 & 1& .& .& 7& 4& 1& 1& 2& 5& 5& 2& 2& 2& 3& .& 2& 2& 1& 1& 1 \\ 
54& 5\,878\,125 & 2& 1& 1& 9& 5& 1& 1& 3& 6& 6& 2& 2& 1& 4& .& 3& 3& 2& 2& 1 \\
 \hline
 \end{array}
 \]}
 \caption{The decomposition matrix of the principal $3$-block $B_0$ 
          of $\mathsf{HN}$.}
 \label{tab:decB0mod3}
\end{table}
\end{landscape}

\subsection{The proof for the block~$B_1$, part 2}
\label{subsec:B1mod3part2}
At the end of Section~\ref{subsec:B1mod3part1} we were left with
44 possible decomposition matrices. Having determined all but one
irreducible Brauer character of $B_1$, it is now elementary that
$a=1=b$, and thus only ten cases remain.

To obtain more information on the last unknown Brauer character of
$B_1$, we condense the tensor product $3\,344\otimes 9\,139$, which
may be done over the field with three elements, making the computation
more efficient. As we have complete knowledge of all other blocks, it
is trivial to only consider the part of the tensor product lying in
$B_1$. It is furthermore sufficient to only condense with respect to the
condensation subgroup $K_2$. We obtain the composition factors
\[
8\,910,\: 16\,929,\: 40\,338,\: 2\times 270\,864,\: 2\times 1\,159\,191,\:
1\,305\,072,\: c1\,047,\]
where $c1\,047$ may possibly be just a composition factor of the 
condensed last simple module restricted to $\Ccal_2$. 
A comparison with the basic set of $B_1$ (see the last column of
Table~\ref{tab:decBSmod3}) yields the relation
\[ \chi_{49}' = 40\,338 + 1\,159\,191 + c1\,047,\]
from which we obtain a Brauer atom of degree $3\,362\,391$. The
latter is a lower bound for the degree of the last irreducible Brauer
character. Inspecting its possible degrees given by the ten possible
decomposition matrices, we see that the obtained lower bound is in fact
maximal, that is we have $c=4$ and $\Phi_3=\Psi'_3-\Phi_7$. Hence the
atom is the character sought, and we give the decomposition matrix of
$B_1$ in Table~\ref{DecMatB1}.
\begin{table}
\[
\begin{array}{r|r|rrrrrrr} \hline
\chi & \chi(1) 
& \Phi_1 & \Phi_2 & \Phi_3 & \Phi_4 & \Phi_5 & \Phi_6 & \Phi_7 \\ \hline
{\bf 8} &   8\,910 & 1 &  . &  . &  . &  . &  . & . \rule[ 0pt]{0pt}{ 13pt} \\
{\bf 10} &  16\,929 &  . &  1 &  . &  . &  . &  . & . \\
{\bf 19} & 270\,864 &  . &  . &  1 &  . &  . &  . & . \\
{\bf 32} &1\,185\,030 &  1 &  1 &  . &  1 &  . &  . & . \\
{\bf 33} &1\,354\,320 &  1 &  . &  . &  . &  1 &  1 & . \\
{\bf 37} &1\,575\,936 &  . &  . &  1 &  . &  1 &  . & . \\
43&2\,784\,375 &  1 &  . &  1 &  1 &  1 &  1 & . \\
{\bf 49} &4\,561\,920 &  . &  . &  . &  1 &  . &  1 & 1 \\
50&4\,809\,375 &  . &  1 &  1 &  1 &  . &  . & 1
\rule[- 7pt]{0pt}{ 5pt} \\ \hline
\end{array}
\]
\caption{The decomposition matrix of the 
$3$-block~$B_1$ of $\mathsf{HN}$}
\label{DecMatB1}
\end{table}

\subsection{The Automorphism Group in Characteristic 3}
The characters of the automorphism group $G.2$ of $G$ fall into
seventeen $3$-blocks, five of which have positive defect: the principal
block $B_0$ of defect $d=6$, having $k=42$ ordinary and $l=22$ modular
characters, two blocks of defect $d=2$, each having $k=9$ ordinary and
$l=7$ modular characters, covering the block $B_1$ of defect $2$ of $G$,
and two blocks of defect $d=1$, each having $k=3$ ordinary and $l=2$
modular characters, covering the block $B_2$ of defect $1$ of $G$.

As $G$ is a normal subgroup of $G.2$, by Clifford's theorem the
irreducible Brauer characters of $G.2$ come in two flavours: if an
irreducible Brauer character of $G$ is not invariant under the action of
the outer automorphism, then induction to $G.2$ yields an irreducible
Brauer character; if the character is invariant, then it possesses two
extensions to $G.2$, where one extension may be derived from the other
by changing the signs of its values on the outer classes. As it
turns out, we may distinguish both extensions by the sign of their value
on the outer class $2C$. Hence, denoting a Brauer character by its
degree, we additionally affix a superscript $+$ or $-$ to identify the
extension.

As each of the blocks $B_1$ and $B_2$ of $G$ is covered by two blocks of
the automorphism group $G.2$ of the same defect, restriction defines a
Morita equivalence between such a block of $G.2$ and $B_1$ respectively
$B_2$ of $G$. Therefore the decomposition matrix for either block of
defect 2 is given in Table~\ref{DecMatB1}, and the decomposition matrix
of either block of defect 1 is given in Table~\ref{DecMatB2}. We may
therefore focus on the principal block $B_0$ of $G.2$.

Owing to the abundance of pairs of characters conjugate under the outer
automorphism, most Brauer characters of $B_0$ are immediately obtained
through induction. We must therefore only determine the extensions of
the invariant Brauer characters.

>From Section~\ref{sec:mod3} it is immediate that the invariant Brauer
characters $1$, $760$, and $3\,344$ of $G$ are liftable. Furthermore, we
may deduce that the restriction of the ordinary character $9\,405^{\pm}$
has the constituents $266$ and $9\,139^{\pm}$. Therefore we only have
to establish the extensions of the four Brauer characters $137\,236$,
$147\,061$, $255\,037$, and $783\,696$ of $G$. Of course, one extension
determines the other, and therefore it suffices to just determine one
member of every pair.

\begin{table}
 \[
 \begin{array}{r|r|r|rrrr}
  \hline & K_1 & K_2 &
  \begin{turn}{90}
   $760^+ \otimes 760^+$
  \end{turn} &
  \begin{turn}{90}
   $760^+ \otimes3\,344^+$
  \end{turn} &
  \begin{turn}{90}
   $1_{4.\mathsf{HS}.2}^{G.2}$
  \end{turn} &
  \begin{turn}{90}
   $3\,344^+ \otimes 3\,344^+$
  \end{turn} \\ \hline
        1^+ &    k1'_1 &   c1'_1 &   1&   .&   4&   5 \\ 
        1^- &    k1'_2 &   c1'_2 &   .&   1&   .&   1 \\ 
      266 &    - &   c6' &   1&   .&   2&   4 \\ 
      760^+ &    k1'_3 &   - &   .&   2&   1&   5 \\ 
      760^- &    k1'_4 &   - &   .&   .&   2&  11 \\ 
     3\,344^+ &    k9'_1 &   - &   2&   1&   6&  11 \\ 
     3\,344^- &    k9'_2 &   - &   .&   3&   .&   4 \\ 
     9\,139^+ &    k8'_1 &   c7'_1 &   2&   .&   5&   8 \\ 
     9\,139^- &    k8'_2 &   c7'_2 &   .&   1&   2&   5 \\ 
    17\,556 &    - &  c32' &   1&   .&   .&   4 \\ 
    24\,528 &    - &  c16'_1 &   .&   .&   2&  10 \\ 
    63\,536 &   k52' &  c16'_2 &   .&   2&   3&   9 \\ 
   137\,236^+ &   k50'_1 &  c40'_1 &   .&   .&   2&   4 \\ 
   137\,236^- &   k50'_2 &  c40'_2 &   .&   1&   .&   . \\ 
   147\,061^+ &   k89' &  c33' &   .&   .&   1&   2 \\ 
   147\,061^- &    &  &   .&   .&   .&   . \\ 
   255\,037^+ &  k173'_1 &  c89'_1 &   1&   .&   3&   5 \\ 
   255\,037^- &  k173'_2 &  c89'_2 &   .&   .&   .&   1 \\ 
   679\,404 &  k276' & c192' &   .&   .&   .&   2 \\ 
   783\,696^+ &   &  &   .&   .&   .&   . \\ 
   783\,696^- &  k322' & c272' &   .&   .&   .&   1 \\ 
   993\,848 &   &  &   .&   .&   .&   . \\ 
   \hline
 \end{array}
 \]
 \caption{Condensation results for $p=3$ for $\mathsf{HN}.2$}
 \label{tab:condresHN2}
\end{table}

\medskip
This is achieved by using condensation again. For the automorphism
group, we choose the same setup as for the simple group, which allows
us to use the computational infrastructure already in place. In other
words, owing to our previous work in Section~\ref{subsec:B0mod3},
we can again use the six elements of \eqref{eq:sixgens} in $G$
plus one element $\alpha$ of the class $2C$ of $G.2$, which is a member of
a set of standard generators of $G.2$ (see \cite{web:ATLAS}). 
Let $\Dcal_i\leq
e_iFGe_i$ be the algebra generated by these elements, hence we have
$\Ccal_i\leq\Dcal_i$. As our condensation subgroups are contained in
$G$ we have $Ve_i\smash\downarrow_{e_iFGe_i} = (V\smash\downarrow_G)
e_i\smash\downarrow_{e_iFGe_i}$ for any $F[G.2]$-module $V$. Hence, the
dimensions of condensed modules are clear from the onset.

Moreover, if a simple $FG$-module $S$ is not invariant, that is
there is a simple $FG$-module $T\not\cong S$ being the image of $S$
under the outer automorphism, and $S\smash\uparrow^{G.2}$ is simple,
then we have $((S\smash\uparrow^{G.2})e_i\smash\downarrow_{\Dcal_i})
\smash\downarrow_{\Ccal_i} =Se_i\smash\downarrow_{\Ccal_i}\oplus
Te_i\smash\downarrow_{\Ccal_i}$, where
$Se_i\smash\downarrow_{\Ccal_i}\not\cong
Te_i\smash\downarrow_{\Ccal_i}$. Thus
$(S\smash\uparrow^{G.2})e_i\smash\downarrow_{\Dcal_i}$ either is simple,
of dimension twice the dimension of $Se_i\smash\downarrow_{\Ccal_i}$,
or has two non-isomorphic constituents extending
$Se_i\smash\downarrow_{\Ccal_i}$ and $Te_i\smash\downarrow_{\Ccal_i}$,
respectively.
If $S$ is invariant and it 
extends to $G.2$ as $S^{\pm}$, then
$(S^{\pm}e_i\smash\downarrow_{\Dcal_i})\smash\downarrow_{\Ccal_i}
=Se_i\smash\downarrow_{\Ccal_i}$ 
implies that both $S^{\pm}e_i\smash\downarrow_{\Dcal_i}$ are simple.
Moreover, since the actions of $\alpha$ on $S^{\pm}$ differ just by sign,
we conclude that the actions of $e_i\alpha e_i$ on $S^{\pm}e_i$ also
differ by sign only. Hence $e_i\alpha e_i$ acts by the zero map on $S^+e_i$
if and only if it does so on $S^-e_i$, and in this case we have
$S^+e_i\smash\downarrow_{\Dcal_i}\cong S^-e_i\smash\downarrow_{\Dcal_i}$.
Otherwise we have
$S^+e_i\smash\downarrow_{\Dcal_i}\not\cong S^-e_i\smash\downarrow_{\Dcal_i}$,
as by Schur's Lemma any isomorphism is scalar and hence commutes with the
action of $e_i \alpha e_i$.
Thus to decide which case occurs
we only have to consider the action of $e_i\alpha e_i$ 
on either of $S^{\pm}e_i$.

The modules considered are given in Table~\ref{tab:condresHN2},
all of whose restrictions to $G$ have already been considered in 
Table~\ref{tab:condresmod3}, where in Table~\ref{tab:condresHN2} 
we only give the composition factors lying in the principal block,
since restriction to $B_0$ of the condensation results is easily achieved.
Thus the matching between the simple modules of the condensed algebras 
$\Dcal_1$ and $\Dcal_2$ can be derived from the matching obtained 
in Section~\ref{subsec:B0mod3}. We moreover see that
$(S\smash\uparrow^{G.2})e_i\smash\downarrow_{\Dcal_i}$ is simple
for all non-invariant simple $FG$-modules occurring, and
$S^{\pm}e_i\smash\downarrow_{\Dcal_i}$ is simple for all invariant simple 
$FG$-modules occurring. Furthermore, we can conclude from the condensation 
results that $S^{+}e_i\smash\downarrow_{\Dcal_i}$ and
$S^{-}e_i\smash\downarrow_{\Dcal_i}$ are non-isomorphic.
The matching of a condensed module to a member of a pair of extensions
follows from the characters of the considered modules
and through the subsequent analysis of the columns of
Table~\ref{tab:condresHN2} as follows. 

The representations needed to condense the modules of
Table~\ref{tab:condresHN2} are again obtained by restricting the
$4\,371$-dimensional representation of the Baby Monster group $\mathsf{B}$
to $G.2$:
\[ 4\,371\smash\downarrow_{G.2} = 1^- + 266 + 760^- + 3\,344^-,\]
from which we immediately obtain $760^+$ and $3\,344^+$.

\medskip
$\blacktriangleright$
To identify the condensed modules belonging to
$760^+$ and $3\,344^+$, we condense them individually to reveal that 
$$ 760^+e_1\mapsto k1'_3,\: 3\,344^+e_1\mapsto k9'_1 ;$$
to distinguish the composition factors
found here from earlier ones notationally we add a dash.
 
\medskip
$\blacktriangleright$ From the character of $760^+ \otimes 760^+$ 
it is clear that the extension $255\,037^+$ is a
constituent: The restriction to $B_0$ of $760^+ \otimes 760^+$
is the sum $\chi_1'+\chi_6'+\chi_8'+\chi_{11}'+\chi_{20}'$. By the
previous section and the above
we know that $\chi_{20}'$ has as its constituents
$3\,344^+$ and extensions of $255\,037$, and $9\,139$. As $\chi_{20}'$ takes
the value $384$ on the class $2C$, and $3\,344^+$ and $9\,139^\pm$ take the
values $36$ and $\pm 51$, respectively, the consituent of $\chi_{20}'$ of
degree $255\,037$ also assumes a positive value on $2C$, i.e.\ it is
$255\,037^+$.

Hence for the component in $B_0$, 
Table~\ref{tab:condresHN2} leaves the two possibilities
\[ 1^+ + 266 + 2\times 3\,344^+ + 17\,556 + 2 \times 9\,139^\pm
+ 255\,037^+,\]
i.e., either $9\,139^+$ or $9\,139^-$ is a constituent. 
Assuming $9\,139^-$ is, it turns out that the character
$\chi_{20}'$ cannot be written as a linear combination of the known Brauer
characters with non-negative integer coefficients. Hence, $9\,139^+$ is a
constituent, from which we obtain $255\,037^+$ explicitly.
Furthermore, we identify the condensed modules
$$ 9\,139^+e_1\mapsto k8'_1,\: 255\,037^+e_1\mapsto k173'_1 .$$

\medskip
$\blacktriangleright$ 
Each of the three modules
$760^+ \otimes 3\,344^+$, $1_{4.\mathsf{HS}.2}^{G.2}$
and $3\,344^+\otimes 3\,344^+$ 
only introduces one new composition factor
which is an extension of a simple $FG$-module.
Hence we may easily deduce their corresponding Brauer characters
$137\,236^-$, $147\,061^+$, and $783\,969^-$, 
where 
\[ 137\,236^-e_1\mapsto k50'_2,\: 147\,061^+e_1\mapsto k89',\: 
783\,969^-e_1\mapsto k322'.\]

The decomposition matrix for the principal $3$-block of $G.2$ is given in
Table~\ref{tab:decB0HN2mod3}.
\begin{landscape}
\begin{table}
\vspace*{-13.7mm} {\tiny
 \[
 \begin{array}{r|r|rrrrrrrrrrrrrrrrrrrrrr}
  \hline
\chi & \chi(1) & \Phi_1 & \Phi_2 & \Phi_3 & \Phi_4 & \Phi_5 & \Phi_6 &
\Phi_7 & \Phi_8 & \Phi_9 & \Phi_{10} & \Phi_{11} & \Phi_{12} &
\Phi_{13} & \Phi_{14} & \Phi_{15} & \Phi_{16} & \Phi_{17} &
\Phi_{18} & \Phi_{19} & \Phi_{20} &\Phi_{21} &\Phi_{22}\\
\hline
  1 &      1 & 1& .& .& .& .& .& .& .& .& .& .& .& .& .& .& .& .& .& .& .& .& . \\ 
  2 &      1 & .& 1& .& .& .& .& .& .& .& .& .& .& .& .& .& .& .& .& .& .& .& . \\ 
  3 &    266 & .& .& 1& .& .& .& .& .& .& .& .& .& .& .& .& .& .& .& .& .& .& . \\ 
  4 &    760 & .& .& .& 1& .& .& .& .& .& .& .& .& .& .& .& .& .& .& .& .& .& . \\ 
  5 &    760 & .& .& .& .& 1& .& .& .& .& .& .& .& .& .& .& .& .& .& .& .& .& . \\ 
  6 &   3\,344 & .& .& .& .& .& 1& .& .& .& .& .& .& .& .& .& .& .& .& .& .& .& . \\ 
  7 &   3\,344 & .& .& .& .& .& .& 1& .& .& .& .& .& .& .& .& .& .& .& .& .& .& . \\ 
  8 &  17\,556 & .& .& .& .& .& .& .& .& .& 1& .& .& .& .& .& .& .& .& .& .& .& . \\ 
 11 &   9\,405 & .& .& 1& .& .& .& .& 1& .& .& .& .& .& .& .& .& .& .& .& .& .& . \\ 
 12 &   9\,405 & .& .& 1& .& .& .& .& .& 1& .& .& .& .& .& .& .& .& .& .& .& .& . \\ 
 15 &  70\,224 & .& .& .& .& .& 1& 1& .& .& .& .& 1& .& .& .& .& .& .& .& .& .& . \\ 
 16 & 131\,670 & 1& 1& .& 1& 1& .& .& .& .& 1& 2& 1& .& .& .& .& .& .& .& .& .& . \\ 
 18 & 214\,016 & 1& .& .& .& 1& 1& .& 1& .& .& .& 1& 1& .& .& .& .& .& .& .& .& . \\ 
 19 & 214\,016 & .& 1& .& 1& .& .& 1& .& 1& .& .& 1& .& 1& .& .& .& .& .& .& .& . \\ 
 20 & 267\,520 & .& .& .& .& .& 1& .& 1& .& .& .& .& .& .& .& .& 1& .& .& .& .& . \\ 
 21 & 267\,520 & .& .& .& .& .& .& 1& .& 1& .& .& .& .& .& .& .& .& 1& .& .& .& . \\ 
 24 & 365\,750 & 1& .& 1& 1& .& 1& .& 1& 1& .& 1& 1& .& .& .& .& 1& .& .& .& .& . \\ 
 25 & 365\,750 & .& 1& 1& .& 1& .& 1& 1& 1& .& 1& 1& .& .& .& .& .& 1& .& .& .& . \\ 
 26 & 749\,056 & .& .& .& 2& 2& .& .& .& .& 1& 2& .& .& .& .& .& .& .& 1& .& .& . \\ 
 29 & 653\,125 & 1& .& .& .& 1& 2& .& 1& 1& .& 1& 1& 1& .& 1& .& 1& .& .& .& .& . \\ 
 30 & 653\,125 & .& 1& .& 1& .& .& 2& 1& 1& .& 1& 1& .& 1& .& 1& .& 1& .& .& .& . \\ 
 31 &1\,312\,500 & 1& 1& .& .& .& .& .& .& .& .& 1& .& .& .& 1& 1& .& .& .& .& .& 1 \\ 
 32 &1\,436\,400 & .& .& 1& 2& 2& 2& 2& 2& 2& 1& 2& 2& .& .& .& .& 1& 1& 1& .& .& . \\ 
 33 &1\,053\,360 & .& .& .& 3& .& .& .& .& .& 1& 2& 1& .& 1& .& .& .& .& .& 1& .& . \\ 
 34 &1\,053\,360 & .& .& .& .& 3& .& .& .& .& 1& 2& 1& 1& .& .& .& .& .& .& .& 1& . \\ 
 40 &1\,361\,920 & .& .& .& 1& 2& 2& 1& 1& 1& .& 2& 1& 1& .& 1& .& 1& .& 1& .& .& . \\ 
 41 &1\,361\,920 & .& .& .& 2& 1& 1& 2& 1& 1& .& 2& 1& .& 1& .& 1& .& 1& 1& .& .& . \\ 
 42 &2\,723\,840 & .& .& .& 2& 2& 1& 1& 1& 1& .& 2& .& .& .& 1& 1& .& .& 2& .& .& 1 \\ 
 49 &2\,375\,000 & 1& 1& 1& 2& 1& 2& 1& 2& 1& .& 2& 1& .& .& 1& 1& 1& .& 1& .& .& 1 \\ 
 50 &2\,375\,000 & 1& 1& 1& 1& 2& 1& 2& 1& 2& .& 2& 1& .& .& 1& 1& .& 1& 1& .& .& 1 \\ 
 51 &2\,407\,680 & 1& 1& 1& .& 1& 1& .& .& 2& .& 2& 2& 1& .& 2& .& .& .& .& 1& .& 1 \\ 
 52 &2\,407\,680 & 1& 1& 1& 1& .& .& 1& 2& .& .& 2& 2& .& 1& .& 2& .& .& .& .& 1& 1 \\ 
 53 &2\,661\,120 & 1& 1& 1& 2& 2& 1& 1& 2& 2& .& 3& 1& .& .& 1& 1& 1& 1& 1& .& .& 1 \\ 
 54 &2\,661\,120 & 1& 1& 1& 2& 2& 1& 1& 2& 2& .& 3& 1& .& .& 1& 1& 1& 1& 1& .& .& 1 \\ 
 59 &3\,200\,000 & .& 1& 1& 2& 1& 1& 1& .& 2& .& 3& 1& .& 1& 2& 1& .& .& 1& 1& .& 1 \\ 
 60 &3\,200\,000 & 1& .& 1& 1& 2& 1& 1& 2& .& .& 3& 1& 1& .& 1& 2& .& .& 1& .& 1& 1 \\ 
 61 &3\,424\,256 & 1& 1& 1& 2& 3& 2& .& 1& 2& 1& 4& 2& 1& .& 2& .& 1& .& 1& 1& .& 1 \\ 
 62 &3\,424\,256 & 1& 1& 1& 3& 2& .& 2& 2& 1& 1& 4& 2& .& 1& .& 2& .& 1& 1& .& 1& 1 \\ 
 65 &4\,156\,250 & .& 1& .& 4& 3& 2& 2& 1& 1& 1& 5& 2& 1& 1& 2& 1& .& .& 2& 1& .& 1 \\ 
 66 &4\,156\,250 & 1& .& .& 3& 4& 2& 2& 1& 1& 1& 5& 2& 1& 1& 1& 2& .& .& 2& .& 1& 1 \\ 
 77 &5\,878\,125 & 1& 1& 1& 4& 5& 3& 2& 1& 2& 1& 6& 2& 1& .& 3& 1& .& .& 3& 1& .& 2 \\ 
 78 &5\,878\,125 & 1& 1& 1& 5& 4& 2& 3& 2& 1& 1& 6& 2& .& 1& 1& 3& .& .& 3& .& 1& 2 \\
 \hline
 \end{array}
 \]}
 \caption{The decomposition matrix of the principal 
$3$-block $B_0$ of $\mathsf{HN}.2$}
 \label{tab:decB0HN2mod3}
\end{table}
\end{landscape}

\bibliographystyle{amsplain}
\bibliography{hnmod2mod3}
\end{document}